\documentclass[twocolumn]{autart}    

\usepackage{graphicx}
\usepackage{lineno,hyperref}
\usepackage{natbib}
\usepackage{tabularx} 
\usepackage{booktabs}

\usepackage{amsmath,amssymb,amsfonts}
\usepackage{graphicx}
\usepackage{amssymb}
\usepackage{eqnarray} 
\usepackage{amsmath}
\usepackage{mathalfa} 
\usepackage[dvips]{epsfig}    
\usepackage{algorithm}
\usepackage{algpseudocode}
\usepackage{psfrag}
\usepackage{amsmath}
\usepackage{stfloats}
\usepackage{amssymb}
\usepackage{dsfont}
\usepackage{xcolor}
\usepackage{bbm}
\usepackage{float}

\DeclareMathOperator{\Tr}{Tr}

\newcommand{\sC}{\mathcal{C}}
\newcommand{\sD}{\mathcal{D}}

\newcommand{\sL}{\mathcal{L}}

\newcommand{\sT}{\mathcal{T}}

\newcommand{\R}{\mathbb{R}}
\newcommand{\E}{\mathbb{E}}
\newcommand{\Z}{\mathbb{Z}}
\newcommand{\K}{\mathbb{K}}

\newcommand{\bv}[1]{\mathbf{#1}}

\newcommand{\plant}[1]{p_{#1}^{(x)}}
\newcommand{\control}[1]{p_{#1}^{(u)}}
\newcommand{\candidatecontrol}[1]{\tilde{p}_{#1}^{(u)}}
\newcommand{\optimalcontrol}[1]{{p_{#1}^{(u)}}^\star}
\newcommand{\modifiedcontrol}[1]{\bar{q}_{#1}^{(u)}}
\newcommand{\data}{\boldsymbol{\Delta}}

\newcommand{\idealplant}[1]{q_{#1}^{(x)}}
\newcommand{\idealcontrol}[1]{q_{#1}^{(u)}}

\newcommand{\KL}{\text{KL}}

\newcommand{\DKL}{{D}_{\KL}}

\newtheorem{Remark}{Remark}
\newtheorem{Theorem}{Theorem}

\newtheorem{Lemma}{Lemma}
\newtheorem{Corollary}{Corollary}
\newtheorem{Assumption}{Assumption}

\newtheorem{Problem}{Problem}

\begin{document}

\begin{frontmatter}

\title{On Convex Data-Driven Inverse Optimal Control for Nonlinear,  Non-stationary and Stochastic Systems\thanksref{footnoteinfo}} 

\author[UniSa]{Emiland Garrabe}\ead{egarrabe@unisa.it},
\author[UniSannio]{Hozefa Jesawada}\ead{jesawada@unisannio.it},
\author[UniSannio]{Carmen Del Vecchio}\ead{c.delvecchio@unisannio.it},
\author[UniSa]{Giovanni Russo}\ead{giovarusso@unisa.it}

\thanks[footnoteinfo]{An early version of this paper with a sketch of the proof for one of the results can be found in \cite{garrabe2023inverse}. EG/HJ are joint first authors. Corresponding author: GR.}

\address[UniSa]{Department of Information and Electrical Engineering and Applied Mathematics, University of Salerno, Italy}
\address[UniSannio]{Department of Engineering, University of Sannio, Benevento, Italy}

\begin{abstract}
This paper is concerned with a finite-horizon inverse control problem,  which has the goal of reconstructing, from observations, the possibly non-convex and non-stationary cost driving the actions of an agent.  In this context, we present a result enabling cost reconstruction by solving an optimization problem that is convex even when the agent cost is not and when the underlying dynamics is nonlinear, non-stationary and stochastic.  To obtain this result, we also study a finite-horizon forward control problem that has randomized policies as decision variables.  We turn our findings into algorithmic procedures and show the effectiveness of our approach via  in-silico and hardware validations.  All  experiments confirm the effectiveness of our approach.
\end{abstract}

\end{frontmatter}

\section{Introduction}

Inverse optimal control/reinforcement learning (IOC/ IRL) refer to the problem of reconstructing the cost driving the actions of an agent from input/output observations \cite{ABAZAR2020119}.  
Tackling IOC/IRL problems is crucial to many scientific domains from engineering, psychology, economics, management and computer science.  
However, a key challenge is that the underlying optimization can become ill-posed even for a linear deterministic dynamics with convex cost.  
Motivated by this,  in a finite-horizon setting, we propose an algorithm enabling cost reconstruction via convex optimization. 
The underlying optimization problem is convex even when the agent cost is not and when the underlying dynamics is nonlinear, non-stationary and stochastic. 
The approach leverages probabilistic descriptions that can be obtained directly from data and/or from first-principles. 
We briefly survey some works related to the results/framework of this paper. For a comprehensive review on IOC/IRL  we refer to  \cite{ABAZAR2020119}.

\noindent{\bf Related Works.} IOC  methods were originally developed to determine control {\em functions} that generate observed outputs \cite{506395}.  More recently,  also driven by the advancements in computational power and the increase in the availability of data, there has been a renewed interest in IOC (and IRL) methods. For stationary Markov Decision Processes ({MDPs}) an approach based on maximum entropy is proposed in \cite{10.5555/1620270.1620297}, with the resulting algorithm based on a backward/forward pass scheme.  {The framework of maximum causal entropy is extended to the infinite time horizon setting in \cite{8115277},  where two formulations are proposed together with convex approximations.} Additionally, following this research stream, \cite{Maximum-Entropy-Multi-Agent} considers linear multi-agent games, while \cite{Levine_Locally_Optimal} obtains results based on the local approximation of the cost/reward. Furthermore, \cite{NIPS2011_c51ce410} also proposes a complementary approach based on the use of Gaussian Processes to model stochastic dynamics.  The resulting algorithm requires matrix inversions and relies on optimization problems that are not  convex in general.  For stationary MDPS, \cite{Finn_Deep_Learning} builds on maximum-entropy and uses deep neural networks to estimate the cost. The approach is benchmarked on manipulation tasks. Also for manipulation tasks, \cite{6630743} considers deterministic nonlinear systems, utilizing path integrals to learn the cost.  In the context of deterministic systems, \cite{SELF2022110242} introduces a model-based IRL algorithm, while \cite{LIAN2022110524} tackles the IRL problem for multiplayer non-cooperative games.  Linearly solvable MDPs are exploited in \cite{IRL-Todorov} and within this framework one can obtain a convex optimization problem to estimate the cost. However, this approach assumes that the agent can specify  the state transition. A risk-sensitive IRL algorithm is proposed in  \cite{Risk-Sensitive} for cost estimates in stationary MDPs assuming that the expert policy belongs to the exponential distribution.  {For LTI deterministic systems,  in \cite{9409781} an IRL algorithm is presented that learns the parameters of a quadratic cost for tracking problems. } In the context of IOC, \cite{NAKANO2023110831} considers stochastic dynamics and proposes an approach to estimate the parameters of a control regularizer for finite-horizon problems.  {For finite-horizon discrete-time linear quadratic regulators,  in \cite{YU2021109636} the cost parameters  are estimated from noisy measurements. }  Also, \cite{Rodrigues} tackles an IOC problem for known nonlinear deterministic systems with quadratic cost  in the input. In \cite{DO2019539} infinite-horizon IOC problems are considered for inverse optimal stabilization and inverse optimal gain assignment for stochastic nonlinear systems driven by L\'evy processes. In \cite{DENG1997151}, stabilization problems in an infinite-horizon setting are considered and it is shown that for every system with a stochastic control Lyapunov function, one can construct a controller which is optimal with respect to some cost.  In \cite{9429728} the cost design problem  is considered and  combinations of input/states compatible with a given value function are studied. Finally,   the approach we  propose relies on a technical result to find the optimal solution for a related finite-horizon forward control problem.  This problem has randomized policies as decision variables and {and it is related to the problems considered in \cite{KARNY19961719,Todorov_pnas,Karny_M+Guy_TV_Sys&Ctr_lett_2006} and references therein.  See} \cite{GARRABE202281} for a survey on this class of sequential decision-making problems.

\subsection{Contributions}\label{sec:contributions}
  We propose an inverse control algorithm, which enables cost reconstruction by solving an optimization problem that is convex even when the agent cost is not convex, non-stationary and when the dynamics is nonlinear, non-stationary and stochastic.  To the best of our knowledge, this is the first algorithm featuring an underlying optimization problem that is guaranteed to be convex in this setting and that,  at the same time: (i) does not require that the agent can specify its state transitions; (ii) does not assume a stationary/discounted setting; (iii) does not require solving forward problems in each iteration and can be used for stochastic systems.  Our key technical contributions are  as follows:
\begin{enumerate}
\item {I}n the finite-horizon setting, we introduce a result to recast the IOC problem into a convex optimization problem.  We show that convexity is guaranteed even when the task cost is not convex and non-stationary. Moreover, the result, which leverages certain probabilistic descriptions  that can be obtained directly from data, does not require assumptions on the underlying dynamics besides the standard Markov assumption;
\item {T}o obtain our result on the inverse problem, we also {study} a related {F}orward {Optimal Control (FOC)} problem.  {For this} finite-horizon problem with randomized policies as decision variables{, we give an expression of the optimal policy, exploited to tackle our IOC problem.  We discuss when this problem has a bounded cost (for the precise definition see Assumption \ref{asn:bounded_cost}) and this also leads to a comprehensive discussion, interesting {\em per se}, relating the problem we consider to other control problem formulations that involve KL Divergence minimization;}
\item {T}he results for the inverse and the forward control problems are turned into algorithmic procedures{.  Leveraging these algorithms, we  illustrate the effectiveness of  our results on problems that involve both continuous and discrete states and actions. Moreover, we compare our IOC result with two state-of-the-art algorithms in the literature. The experiments illustrate that our algorithm outperforms those from the literature achieving a lower discrepancy between real and reconstructed costs  with significantly lower computation times.  Code available at \url{https://tinyurl.com/46uccxsf}};
\item {T}he algorithms are validated via a hardware test-bed (the documented code for the experiments is available at \url{https://tinyurl.com/46uccxsf}). Namely, we use our algorithms to reconstruct the cost of robots navigating in an environment with the goal of reaching a desired destination{. We consider both an obstacle-free environment (Scenario $1$) and an environment populated with obstacles (Scenario $2$).  The experimental results,  which are obtained for both continuous and discrete variables,  demonstrate the effectiveness of our approach.  Moreover,  t}he paper also includes a running example involving the swing-up of a pendulum. The example illustrates some of the key
algorithmic details of our results (code also available at \url{https://tinyurl.com/46uccxsf}).
\end{enumerate}

{While our results are related to works on IRL and IOC that leverage maximum entropy and likelihood arguments\footnote{{Notably, the exponential policies arising from our FOC problem and exploited to tackle the IOC are also {\em assumed} in these works.  This  is further discussed throughout the paper}} such as \cite{10.5555/1620270.1620297,IRL-Todorov,8115277}, this work extends these in several ways. First,  we do not necessarily require that the MDP is stationary as instead assumed in these papers. Second,  we show that the underlying optimization problem that we obtain  to reconstruct the cost is convex even when the dynamics is stochastic and the agent cost is not convex and non-stationary.  In \cite{10.5555/1620270.1620297} convexity is guaranteed only when the MDP is deterministic.  In \cite{8115277} a convex IRL formulation is obtained but this relies on an approximation of the original problem.  Instead, we do not rely on approximations to obtain convexity.  Third,  the approach we propose can be used for both continuous and discrete action/state spaces, while \cite{10.5555/1620270.1620297,IRL-Todorov,8115277} tackle IRL/IOC problems with discrete variables. Moreover, we do not require that the agent is able to specify its state transitions as instead assumed in \cite{IRL-Todorov}.  Additionally,  in \cite{10.5555/1620270.1620297,8115277}, policy evaluation is required at each iteration of the algorithm and this is not required with our approach.  Finally,  we also deploy the policies from the FOC problem on a setup involving real hardware. Consequently,  we also show that the cost can be effectively reconstructed from these policies.  This was not done in the above works.}

\section{Mathematical {Background}}\label{sec:math_preliminaries_formulation}

Sets are in {\em calligraphic} and vectors in {\bf bold}.  We let $\K$ be either $\R$ or $\Z$.  A random variable is denoted by $\mathbf{V}$ and its realization is $\mathbf{v}$. We denote the \textit{probability mass function} (pmf, for discrete variables) or \textit{probability density function} (pdf, for continuous variables) of $\mathbf{V}$  by $p(\mathbf{v})$ and we let $\sD$ be the convex subset of pdfs/pmfs.  In what follows, we simply say that $p(\mathbf{v})$ is a probability function (pf). Whenever we take the sums/integrals involving pfs we always assume that they exist.  The  expectation of a function $\mathbf{h}(\cdot)$ of a discrete $\mathbf{V}$ is denoted $\E_{{p}}[\mathbf{h}(\mathbf{V})]:=\sum_{\bv{v}}\mathbf{h}(\mathbf{v})p(\mathbf{v})$, where the sum is over the support of $p(\mathbf{v})$; whenever it is clear from the context, we  omit the subscript in the sum (for continuous variables the summation is replaced by the integral on the support of $p(\mathbf{v})$). The {joint} pf of $\mathbf{V}_1$ and $\mathbf{V}_2$ is denoted by  $p(\mathbf{v}_1,\mathbf{v}_2)$ and the {conditional} pf of $\mathbf{V}_1$ with respect to (w.r.t.) $\mathbf{V}_2$ is $p\left( \mathbf{v}_1\mid   \mathbf{v}_2 \right)$. Countable sets are denoted by $\lbrace w_k \rbrace_{k_1:k_n}$, where $w_k$ is the generic set element, $k_1$ ($k_n$) is the index of the first (last) element and  $k_1:k_n$ is the set of consecutive integers between (including) $k_1$ and $k_n$. A pf of the form $p(\bv{v}_0,\ldots,\bv{v}_N)$ is compactly written as $p_{0:N}$ (by definition $p_{k:k} := p_k(\bv{v}_k)$).  {We compactly write conditional pfs of the form $p_k(\bv{y}_k\mid\bf{z}_{k-1})$ as $p^{(y)}_{k\mid k-1}$}.
Also,  functionals are denoted by capital calligraphic characters with arguments within curly brackets. We make use of the Kullback-Leibler (KL, \cite{KL_51}) divergence, a measure of the proximity of the pair of pmfs  $p(\mathbf{v})$ and $q(\mathbf{v})$, defined for discrete variables as
${D}_{\KL}\left(p \mid\mid q \right):= \sum_{\bv{v}} p(\bv{v}) \ln\left( {p(\bv{v})}/{q(\bv{v})}\right)$.  {F}or continuous variables the sum is replaced by the integral.  {The ${D}_{\KL}\left(p \mid\mid q \right)$ is finite only if the support of $p(\bv{v})$ is contained in the support of $q(\bv{v})$. That is,  ${D}_{\KL}\left(p \mid\mid q \right)$ is finite only if $p(\bv{v})$ is absolutely continuous with respect to (w.r.t.) $q(\bv{v})$,  see e.g.,  \cite[Chapter $8$]{10.5555/1146355}}. We recall the chain rule for the KL divergence:

\begin{Lemma}\label{lem:splitting_property}	
Let $\bv{V}$ and $\bv{Z}$ be two random variables  and let $f(\bv{v},\bv{z})$ and $g(\bv{v},\bv{z})$ be two joint pfs. Then:
\begin{multline*}
{D}_{\KL}
\left( f(\bv{v},\bv{z}) \mid\mid g(\bv{v},\bv{z})  \right) = 
\\ {D}_{\KL}
\left( f(\bv{v}) \mid\mid g(\bv{v})  \right)
+ \mathbb{E}_{f(\bv{v})}
\left[	
{D}_{\KL} 
\left(f(\bv{z}\mid\bv{v})\mid\mid g(\bv{z}\mid\bv{v}) 
\right)
\right].
\end{multline*}
\end{Lemma}

\subsection{Problems Set-up}\label{sec:set-up}

We let $\bv{X}_k\in\mathcal{X}\subseteq\K^n$ be the system state at time-step $k$ and $\bv{U}_k\in\mathcal{U}\subseteq\K^p$ be the control input at time-step $k$.  The time indexing is chosen so that the system transitions to $\bv{x}_k$ when $\bv{u}_k$ is applied.  That is,  by making the Markov assumption,  the possibly non-stationary, nonlinear stochastic dynamics for the system under control is {given by the pf $p^{(x)}_k	\left(	\mathbf{x}_k \mid \mathbf{u}_k, \mathbf{x}_{k-1} \right)$.  In what follows, we make use of the short-hand notation $\plant{k\mid k-1}$ to denote this pf, which is simply termed} as target pf.  

\begin{Remark}
In the running example we obtain the target pf from data.  We let $\data_k:=(\bv{x}_{k-1},\bv{u}_k)$ be a data pair.  Also, $\data_{0:N}:=(\{\data_k\}_{1:N},\bv{x}_N)$ is the {\em dataset} collected over  $\sT:={0:N}$.  A  {\em database}  is a collection of datasets.
\end{Remark}

To formalize the control problems we introduce the pf:
\begin{equation}\label{eqn:behavior}
\begin{split}
p_{0:N} & = p_{0}\left(\bv{x}_0\right)\prod_{k=1}^N p_{k\mid k-1} = p_{0}\left(\bv{x}_0\right)\prod_{k=1}^N 
\plant{k\mid k-1}
\control{k\mid k-1},
\end{split}
\end{equation}
where $\control{k\mid k-1}:=p^{(u)}_k	\left( \mathbf{u}_k \mid \mathbf{x}_{k-1} \right)$ is a randomized policy and initial conditions are embedded via the prior $ p_{0}\left(\bv{x}_0\right)$.  Also, we use the shorthand notation $p_{k\mid k-1}:=\plant{k\mid k-1}\control{k\mid k-1}= p_{{k}}\left(\bv{x}_k,\bv{u}_k\mid\bv{x}_{k-1}\right)$. 

\begin{Remark}
{At each $k$, the pf $p^{(x)}_k	\left(	\mathbf{x}_k \mid \mathbf{u}_k, \mathbf{x}_{k-1} \right)$ describes the behavior of the system given $\bv{x}_{k-1}$ and $\bv{u}_k$;} $p_{0:N}$ describes in probabilistic terms the evolution of closed-loop system when, at each $k$, a given policy, say  $p^{(u)}_k	\left( \mathbf{u}_k \mid \mathbf{x}_{k-1} \right)$, is used.  With the forward control problem {(}Section \ref{sec:main_problem}{)} we aim to design the policy. This is then exploited to tackle the inverse problem {(}Section \ref{sec:cost_estimation_formalization}{)}.
\end{Remark}

\subsection{The Forward Control Problem}\label{sec:main_problem}
We let $c_k:\mathcal{X}\rightarrow\R$ be the cost, at time-step $k$, associated to  $\bv{x}_k$. Then, the expected cost incurred when the system is in state $\bv{x}_{k-1}$ and input $\bv{u}_k$ is applied is given by $\E_{p^{(x)}_{k\mid k-1}}\left[c_k(\bv{X}_k)\right]$.
The forward control problem considered in this paper is formalized with the following:
\begin{Problem}\label{prob:main}
Given  a joint pf 
\begin{equation*}\label{eqn:q_joint}
    q_{0:N} := q_{0}\left(\bv{x}_0\right)\prod_{k=1}^N q^{(x)}_k	\left(	\mathbf{x}_k \mid \mathbf{u}_k, \mathbf{x}_{k-1} \right)
q^{(u)}_k	\left( \mathbf{u}_k \mid \mathbf{x}_{k-1} \right),
\end{equation*}
find the sequence of pfs, $\left\{\optimalcontrol{k\mid k-1}\right\}_{1:N}$, such that:
\begin{equation}\label{eqn:main_problem}
    \begin{aligned}
{\left\{\optimalcontrol{k\mid k-1}\right\}_{1:N}}\in & \underset{\left\{ \control{k\mid k-1}\right\}_{1:N}}{\text{arg min}}
    \bigg\{\DKL\left(p_{0:N}\mid\mid q_{0:N}\right)\\
    & + \sum_{k=1}^N \E_{\bar{p}_{k-1:k}}\left[\E_{p^{(x)}_{k\mid k-1}}\left[c_k(\bv{X}_k)\right]\right]\bigg\} \\
   &  s.t.  \  \control{k\mid k-1}\in\sD \ \ \forall k\in\sT.
    \end{aligned}
\end{equation}
where  $\bar{p}_{k-1:k}:=p_{k-1}(\bv{x}_{k-1},\bv{u}_k)$. 
\end{Problem}

The solution of Problem \ref{prob:main} is a sequence of randomized policies.  At each $k$, the control input applied to the system, i.e. $\bv{u}^\star_k$,  is sampled from  $\optimalcontrol{k\mid k-1}$.  
\begin{Remark}\label{rem:regularizer}
In Problem \ref{prob:main}, minimizing the first term in the cost functional amounts at minimizing the discrepancy between $p_{0:N}$ and $q_{0:N}$. Hence, this first term can be thought of as a regularizer,  biasing the behaviour of the closed loop system towards $q_{0:N}$.  {In Remark \ref{rem:regularizer_follow-up} we discuss  how the optimal solution of the problem changes when a regularization weight is introduced in the formulation.} 
\end{Remark}
Typically,  $q_{0:N}$ can be a passive dynamics\cite{9029512,Todorov_pnas} or capture  behaviors from demonstrations  \cite{gagliardi2020probabilistic}.  {The pfs for the dynamics in Problem \ref{prob:main} can be derived both directly from the data (see e.g.,  our running example) and from first-principles. This setting, which foresees randomized policies as decision variables, is widely adopted in the literature, see e.g. \cite{9029512,Karny_M+Guy_TV_Sys&Ctr_lett_2006,6716965,Todorov_pnas} and \cite{GARRABE202281} for a survey across learning and control.  Moreover, as we shall see, the class of policies obtained as optimal solution of Problem \ref{prob:main} also naturally arises in the context of IOC/IRL, see e.g., \cite{8115277,Todorov_pnas,10.5555/1620270.1620297}.  For this last point,  see also Remark \ref{rem:exponential_policy}.}
  
For our derivations, it is also useful to define  $\idealplant{k\mid k-1} := q^{(x)}_k	\left(	\mathbf{x}_k \mid \mathbf{u}_k, \mathbf{x}_{k-1} \right)$, $\idealcontrol{k\mid k-1}:=q^{(u)}_k	\left( \mathbf{u}_k \mid \mathbf{x}_{k-1} \right)$ and $q_{k\mid k-1}:=\idealplant{k\mid k-1}\idealcontrol{k\mid k-1} = q_{{k}}\left(\bv{x}_k,\bv{u}_k\mid\bv{x}_{k-1}\right)$.  When emphasizing the cost dependency on the decision variables in Problem \ref{prob:main}, we use the notation $\mathcal{J}\left\{p_{k\mid k-1}^{(u)}\right\}_{1:N}.$
\begin{Assumption}\label{asn:bounded_cost}
For Problem \ref{prob:main}, there exists some  feasible $\left\{\tilde p_{k\mid k-1}^{(u)}\right\}_{1:N}$ such that $\mathcal{J}\left\{\tilde p_{k\mid k-1}^{(u)}\right\}_{1:N}$ is bounded{, that is,  there exist some $\underbar{J}$ and $\bar{J}$ such that $-\infty < \underbar{J} \le \mathcal{J}\left\{\tilde p_{k\mid k-1}^{(u)}\right\}_{1:N}\le \bar{J} < +\infty$.}
\end{Assumption}
{We now discuss Assumption \ref{asn:bounded_cost}. To this aim, we  define
\begin{equation}\label{eqn:twisted_ref}
\tilde q_{0:N} := \frac{q_{0:N}\exp\left(-\sum_{k=1}^N\tilde c_k(\bv{x}_{k-1},\bv{u}_k)\right)}{\sum_{\bv{x}_{0:N},\bv{u}_{1:N}}q_{0:N}\exp\left(-\sum_{k=1}^N\tilde c_k(\bv{x}_{k-1},\bv{u}_k)\right)},
\end{equation}
with $\tilde c_k(\bv{x}_{k-1},\bv{u}_k):= \E_{p^{(x)}_{k\mid k-1}}\left[c_k(\bv{X}_k)\right]$.  By leveraging the Markov assumption, the above pf can be decomposed as
$$
\tilde q_{0:N} := \tilde q_{0}\left(\bv{x}_0\right)\prod_{k=1}^N \tilde{q}^{(x)}_k	\left(	\mathbf{x}_k \mid \mathbf{u}_k, \mathbf{x}_{k-1} \right)
\tilde{q}^{(u)}_k	\left( \mathbf{u}_k \mid \mathbf{x}_{k-1} \right),
$$
and again we use the shorthand notation $\tilde{q}^{(x)}_{k\mid k-1}$ and $\tilde{q}^{(u)}_{k\mid k-1}$ to denote $\tilde{q}^{(x)}_k	\left(	\mathbf{x}_k \mid \mathbf{u}_k, \mathbf{x}_{k-1} \right)$ and $
\tilde{q}^{(u)}_k	\left( \mathbf{u}_k \mid \mathbf{x}_{k-1} \right)$, respectively. Also, we let $\tilde{q}_{k\mid k-1}:=  \tilde{q}^{(x)}_{k\mid k-1}\tilde{q}^{(u)}_{k\mid k-1}$ and note that $\mathcal{J}\left\{\tilde p_{k\mid k-1}^{(u)}\right\}_{1:N}$ can be  written as
\begin{equation}\label{eqn:twisted_cost}
\DKL\left(p_{0:N}\mid\mid \tilde q_{0:N}\right) - \sum_{\bv{x}_{0:N},\bv{u}_{1:N}}q_{0:N}\exp\left(-\sum_{k=1}^N\tilde c_k(\bv{x}_{k-1},\bv{u}_k)\right).
\end{equation}
Then,  we are ready to give the following:
\begin{prop}\label{prop:assumption}
Consider Problem \ref{prob:main} and assume that:
\begin{enumerate}
\item[(i)] there exist non-negative constants $H_k<+\infty$, $k=0:N$ such that $\DKL(p_0(\bv{x}_0)\mid\mid \tilde q_0(\bv{x}_0)) \le H_0$ and, $\forall k=1:N$, $\DKL\left(p^{(x)}_{k\mid k-1}\mid\mid \tilde{q}^{(x)}_{k\mid k-1}\right) \le H_k$ ;
\item[(ii)]  there exist some $\bar e$ such that
 $$\E_{q_{0:N}}\left[\exp\left(-\sum_{k=1}^N\tilde c_k(\bv{x}_{k-1},\bv{u}_k)\right)\right] \le \bar{e} < +\infty.
 $$
\end{enumerate}
Then, Assumption \ref{asn:bounded_cost} is satisfied with $\left\{\tilde p_{k\mid k-1}^{(u)}\right\}_{1:N} = \left\{\tilde q_{k\mid k-1}^{(u)}\right\}_{1:N}$ and with $\underbar{J} = -\bar{e}$ and $\bar{J} = \sum_{k=0}^NH_k$.
\end{prop}
\begin{pf} See the appendix.
\end{pf}
\begin{Remark}
In Proposition \ref{prop:assumption},  condition (i) is satisfied when $p_0(\bv{x}_0)$ and $p^{(x)}_{k\mid k-1}$ are absolutely continuous w.r.t. $\tilde q_0(\bv{x}_0)$ and $\tilde q^{(x)}_{k\mid k-1}$, respectively.  Condition (ii) is satisfied whenever the agent cost $c_k(\cdot)$ is lower bounded.
\end{Remark}
\begin{Remark}\label{rem:KL_min_link}
From the above discussion, the problem in \eqref{prob:main} has the same minimizer (but different optimal cost) as the problem $\min_{\left\{ \control{k\mid k-1}\right\}_{1:N}}\DKL\left(p_{0:N}\mid\mid \tilde q_{0:N}\right)$.  These cost functionals are also considered in the literature on probabilistic design \cite{Karny_M+Guy_TV_Sys&Ctr_lett_2006},  entropically regularized optimal transport \cite{Nutz_2020} and Shroedinger Bridges \cite{doi:10.1137/20M1339982}. 
\end{Remark}}

\subsection{The Inverse Control Problem}\label{sec:cost_estimation_formalization}

The inverse control problem we consider consists in recovering the cost-to-go for the agent, say $\bar{c}_k(\cdot)$, and the agent cost $c_k(\cdot)$ given a set of observed states and inputs sampled from $\plant{k\mid k-1}$ and from the agent policy.  
In what follows, we denote by $\hat{\bv{x}}_k$ and $\hat{\bv{u}}_k$ the observed state and control input at time-step $k$.  We also make the following:
\begin{Assumption}\label{asn:inverse_problem}
The cost-to-go is expressed as a linear combination of features.  That is,  $\bar{c}_k\left({\bv{x}}_k\right) = -{\bv{w}_k}^T\bv{h}\left({\bv{x}}_k \right)$, where $\bv{h}(\bv{x}_k):=[{h}_1(\bv{x}_k),\ldots,{h}_F(\bv{x}_k)]^T$ is the features vector and ${h}_i:\mathcal{X}\rightarrow\R$  are known functions, $i\in 1:{F}$, and $\bv{w}_k:=[w_{k,1},\ldots,w_{k,F}]^T$ is a vector of weights.
\end{Assumption}
The assumption is rather common in the literature see e.g., \cite{10.5555/1620270.1620297,Maximum-Entropy-Multi-Agent,6630743,SELF2022110242,IRL-Todorov}. With our results in Section \ref{sec:learning_cost} we propose an approach based on maximum likelihood to recover $\bar{c}_k(\cdot)$  and $c_k(\cdot)$. See e.g., \cite{9661376} for a maximum likelihood framework for linear systems in the context of data-driven control.

\section*{Running Example: Pendulum Control}

We consider the control of a pendulum and in this first part of the running example we introduce the setting\footnote{{All  results were obtained using a laptop with an 12th Gen Intel Core i5-12500H 2.50 GHz processor and 18 GB of RAM}}. The forward  control problem consists in stabilizing  a pendulum on its unstable equilibrium. The cost used to formulate this forward control problem (tackled in the next part of the example) is: 
\begin{equation}\label{eqn:state_cost}
    {c}(\bv{x}_{k}) = (\theta_{k} - \theta_{d})^{2}+0.01(\omega_{k} - \omega_{d})^{2},
\end{equation}
with $\theta_{d}=0$ and $\omega_{d}=0$ ($\theta_{d}=0$ corresponds to the unstable equilibrium of the pendulum). With the inverse control problem, tackled in the last part of this running example, we recover the cost from observed states and control inputs.  

The pendulum dynamics, which is only used to generate data, is given by:
\begin{equation}\label{Pendulum_Dynamics}
    {\begin{aligned}
        \theta_{k} &= \theta_{k-1} + {\omega}_{k-1} dt + W_{\theta} \\
    {\omega}_{k} &
    = {\omega}_{k-1} + \left( \frac{g}{l}\sin(\theta_{k-1})+\frac{u_{k}}{ml^{2}}\right)dt + W_{{\omega}}{,}
    \end{aligned}}
\end{equation}
where $\theta_{k}$ is the angular position, $\omega_{k}$ is the angular velocity and $u_{k}$ is the torque applied on the hinged end. The parameter $l$ is the rod length, $m$ is the mass of the pendulum, $g$ is the gravity and $dt = 0.1$s is the discretization step. Also, $W_{\theta}$ and $W_{{\omega}}$ are sampled from Gaussians, i.e.  $W_{\theta}\sim \mathcal{N}(0,0.05)$ and $W_{{\omega}}\sim \mathcal{N}(0,0.1)$.  We let $\bv{X}_{k}:= [\theta_{k},\omega_{k}]^T$ and $u_{k}\in\mathcal{U}$ with $\mathcal{U}:=[-2.5,2.5]$.  The parameters were chosen as in  \cite{GARRABE202281} so that the {target} pendulum, i.e. the pendulum that we want to control, had parameters $m=1$kg, $l=0.6$m, while the reference pendulum (from which the pfs $\idealplant{k \mid k-1}$ and $\idealcontrol{k \mid k-1}$ are extracted) had $m = 0.5$kg, $l = 0.5$m. To illustrate the application of our results in both the continuous and the discrete settings, we built both pdfs and pmfs for the target and reference pendulum.  These were built  from data obtained by simulating (\ref{Pendulum_Dynamics}). The process we followed is outlined below (see \url{https://tinyurl.com/46uccxsf} for the details). 

\noindent {\bf Target pendulum.} For the discrete setting, we estimated the empirical pmf,  $\plant{k \mid k-1}$, for the target pendulum. To do so, we set $\mathcal{X}:=\left[-\pi,\pi\right]\times[-5,5]$ and: (i) built a database of $10000$ simulations of $100$ time-steps each (at each step of the simulations the control input was sampled from the uniform distribution); (ii) discretized the set $\mathcal{X}$  in $50 \times 50$ bins; (iii) used the histogram filter to estimate the empirical pmf from the data.  The algorithm for the histogram filter can be found in e.g.,  \cite{GARRABE202281} that also provides the related documented code.  Instead, for the continuous setting we estimated the pdf $\plant{k \mid k-1}$ via Gaussian Processes. To this aim, we: (i) used a database of $30$ simulations and $100$ time-steps to build a prior (again, at each $k$ the control input was sampled from the uniform distribution); (ii) picked the covariance function as a squared exponential kernel. By doing so,  at each $k$, $\plant{k \mid k-1}$ was a Gaussian with mean and variance inferred from the data.

\noindent {\bf Reference pendulum.} 	We estimated both a continuous and a discrete $\idealplant{k \mid k-1}$ and this was done by following the process described above for the target pendulum. Instead,  $\idealcontrol{k \mid k-1}$ was obtained as in \cite{GARRABE202281} by adding Gaussian noise to a Model Predictive Control (MPC) policy (able to stabilize the unstable equilibrium of the reference pendulum) and subsequently discretizing this pf. 

\section{Main Results}\label{sec:theory_calculations} 
We first present a result tackling the forward control problem by giving the optimal solution of Problem \ref{prob:main}. Then, this result is used to address the inverse control problem.  To streamline the presentation,  results are stated for discrete variables. The proofs for continuous variables are omitted here for brevity.

\subsection{Tackling the Forward Control Problem}\label{sec:optimal_policy}
With the next result we give the solution to  Problem \ref{prob:main}.
\begin{Theorem}\label{thm:prob_main}
Consider Problem \ref{prob:main} and let Assumption \ref{asn:bounded_cost} hold. Then:
\begin{itemize}
\item[(i)] the problem has the unique solution $\{{\optimalcontrol{k\mid k-1}}\}_{1:N}$ , with 
\begin{equation}\label{eqn:optimal_solution_statement}
\begin{split}
& \optimalcontrol{k\mid  k-1} =\frac{\bar{p}_{k\mid  k-1}^{(u)}\exp\left(-\E_{\plant{k\mid  k-1}}\left[\bar{c}_{k}(\bv{X}_{k})\right]\right)}{\sum_{\bv{u}_{k}}\bar{p}_{k\mid  k-1}^{(u)}\exp\left(-\E_{\plant{k\mid  k-1}}\left[\bar{c}_{k}(\bv{X}_{k})\right]\right)},
\end{split}
\end{equation}
where
\begin{equation}\label{eqn:p_bar}
\bar{p}_{k\mid  k-1}^{(u)}:=\idealcontrol{k\mid  k-1}\exp\left(-\DKL\left(\plant{k\mid  k-1}\mid  \mid  \idealplant{k\mid  k-1}\right)\right),
\end{equation}
where $\bar{c}_{k}:\mathcal{X}\rightarrow\R$ is {given by} 
\begin{subequations}\label{eqn:backward_recursion}
\begin{equation}\label{eqn:backward_recursion_a}
\begin{split}
 \bar{c}_{k}(\bv{x}_{k}) & = c_{k}(\bv{x}_{k}) -\hat{c}_{k}(\bv{x}_{k}),\\
 \end{split}
\end{equation}
\mbox{{with $\hat{c}_{k}:\mathcal{X}\rightarrow\R$} obtained via the  backward recursion}
\begin{equation}
\begin{split}
 {\hat{c}_{N}(\bv{x}_{N})} & {= 0},\\
\hat{c}_{k}(\bv{x}_{k}) & = \ln\left(\E_{\idealcontrol{k+1\mid k}}\left[\exp\left(-\DKL\left(\plant{k+1\mid k}\mid\mid\idealplant{k+1\mid k}\right)\right.\right.\right.\\
&\left.\left. \left.-\E_{\plant{k+1\mid k}}\left[\bar{c}_{k+1}(\bv{X}_{k+1})\right]\right)\right]\right), \ \ {k = 1:N-1;}\\
\end{split}
\end{equation}
\end{subequations}
\item[(ii)] the  minimum is
$-\sum_{k=1}^N\E_{\bar{p}_{k-1}}\left[\hat{c}_{k-1}(\bv{X}_{k-1})\right],$
where $\bar{p}_{k-1}:=p_{k-1}(\bv{x}_{k-1})$.
\end{itemize}
\end{Theorem}
{\begin{Remark}\label{rem:regularizer_follow-up}
It is often convenient to include in the formulation of the problem in \eqref{eqn:main_problem} a regularization weight, say $\varepsilon >0$, for the KL divergence component of the cost.  In this case,  the problem in \eqref{eqn:main_problem} becomes
\begin{equation}\label{eqn:regularized_problem}
    \begin{aligned}
{\left\{\optimalcontrol{k\mid k-1}\right\}_{1:N}}\in & \underset{\left\{ \control{k\mid k-1}\right\}_{1:N}}{\text{arg min}}
    \bigg\{\varepsilon\DKL\left(p_{0:N}\mid\mid q_{0:N}\right)\\
    & + \sum_{k=1}^N \E_{\bar{p}_{k-1:k}}\left[\E_{p^{(x)}_{k\mid k-1}}\left[c_k(\bv{X}_k)\right]\right]\bigg\} \\
   &  s.t.  \  \control{k\mid k-1}\in\sD \ \ \forall k\in\sT.
    \end{aligned}
\end{equation}
The optimal solution of \eqref{eqn:regularized_problem} is again given by Theorem \ref{thm:prob_main} with the only difference being in \eqref{eqn:backward_recursion_a} which now becomes $\bar{c}_{k}(\bv{x}_{k})= c^{(\varepsilon)}_{k}(\bv{x}_{k}) -\hat{c}_{k}(\bv{x}_{k})$, with $c^{(\varepsilon)}_k(\bv{X}_k):=c_k(\bv{X}_k)/\varepsilon$.
\end{Remark}}
\begin{pf}
The proof,  which is by induction, is {based on the steps of} 
\cite{gagliardi2020probabilistic,9244209}.  {In \cite{gagliardi2020probabilistic} the control problem had no $c_k(\cdot)$ while in \cite{9244209} the cost was included, however the decision variables were not randomized policies.  We give a self-contained proof in the appendix.}
\end{pf}

\begin{Remark}\label{rem:exponential_policy}
The optimal solution  in Theorem \ref{thm:prob_main} has an exponential twisted kernel \cite{6716965} and this structure is exploited to prove our result on the inverse problem.  This class of policies,  also known as soft-max/Boltzmann policies,  are often assumed in IRL/IOC works.  {For example,  the policies in \cite{Risk-Sensitive,10.5555/1620270.1620297,8115277} can be derived from \eqref{eqn:optimal_solution_statement} and \eqref{eqn:p_bar} when $\idealcontrol{k\mid  k-1}$ is  uniform and $\plant{k\mid  k-1}$ is equal to $ \idealplant{k\mid  k-1}$.  The policy from which the agent cost is reconstructed in these works is in fact a soft-max.}
\end{Remark}
\begin{Remark}\label{rem:constraints}
{In Problem \ref{prob:main} the space of the policies is unconstrained.  However,  from \eqref{eqn:optimal_solution_statement} and \eqref{eqn:p_bar} we get that the optimal solution given in Theorem \ref{thm:prob_main}  is zero whenever $q^{(u)}_{k\mid k-1} = 0$.  That is,  $q^{(u)}_{k\mid k-1}$  can be used to {\em enforce} that the optimal solution is zero for specific $\bv{u}_k$'s. In turn,  this can be useful when $\bv{u}_k$ has bound constraints.}
\end{Remark}

\subsubsection{Turning Theorem \ref{thm:prob_main} into an Algorithm}\label{sec:results}
\begin{figure}
    \centering
    \includegraphics[width=0.45\columnwidth]{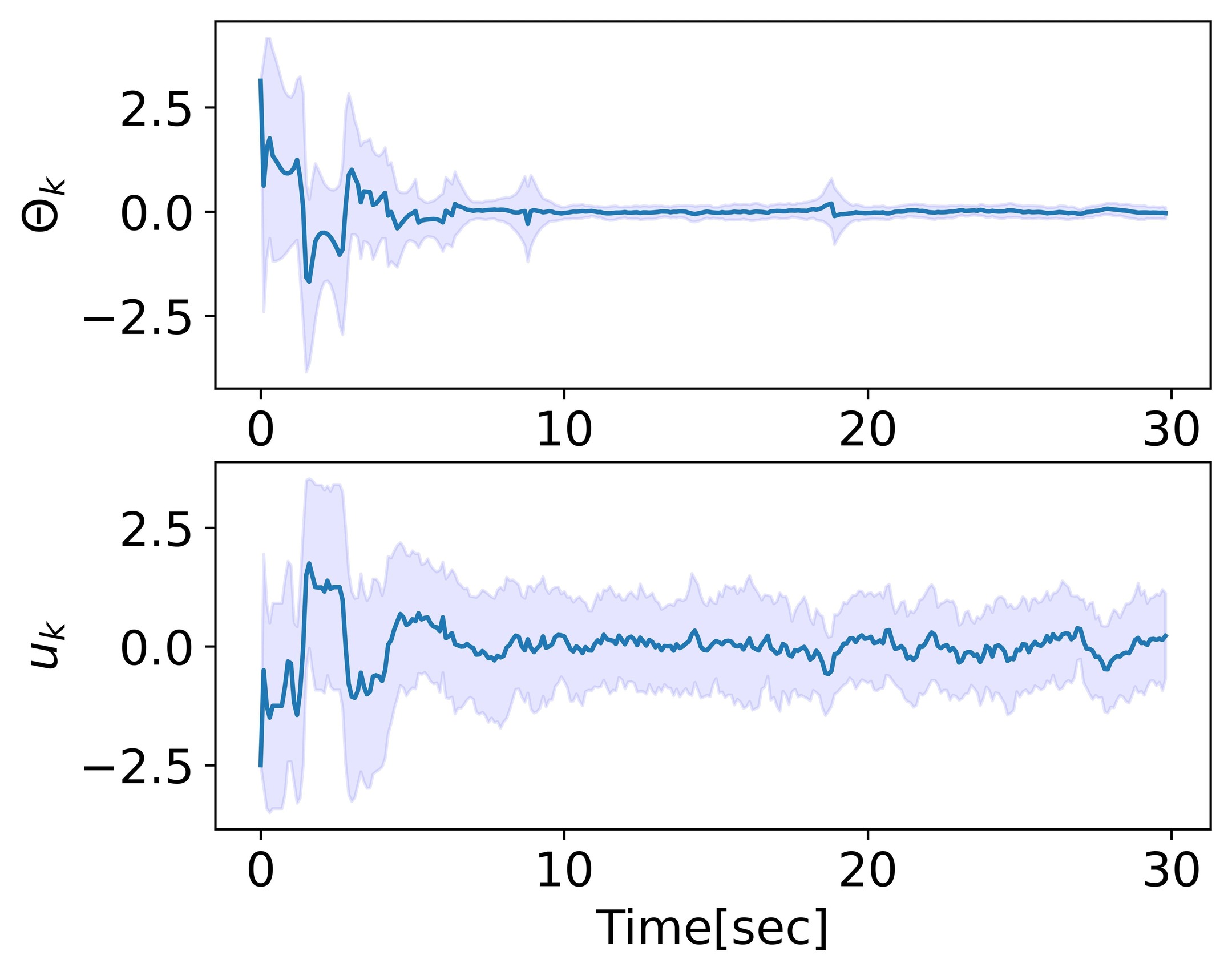}
    \includegraphics[width=0.45\columnwidth]{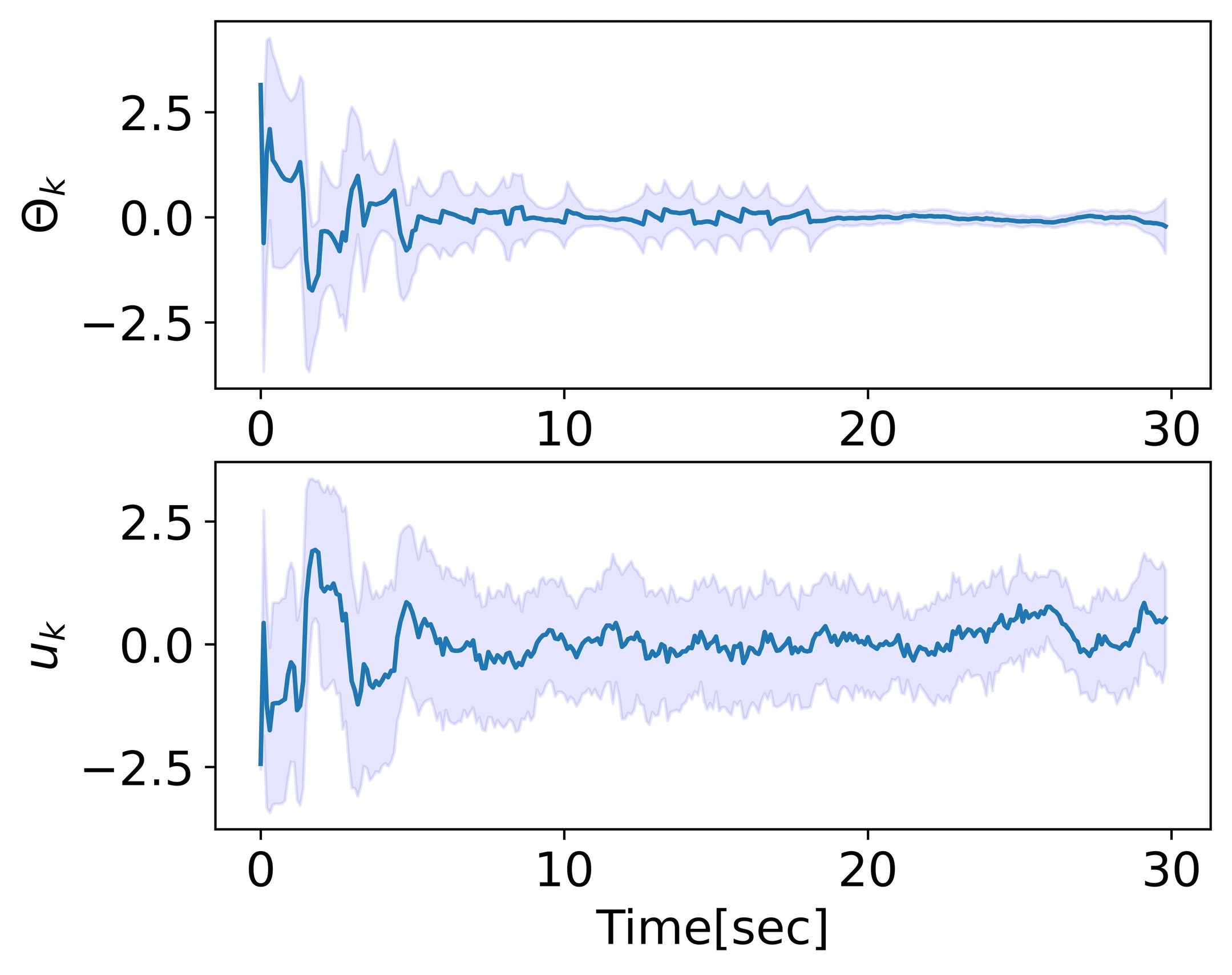}
    \caption{Target pendulum angular position and corresponding control input.  Results obtained when: (i) pfs are discrete,  estimated via the histogram filter (left panels); (ii) pfs are estimated via Gaussian Processes (right panels).  Panels obtained from $20$ simulations. Bold lines represent the mean and the shaded region is confidence interval corresponding to the standard deviation. }
    \label{fig:sim_forward_pendulum}
\end{figure}

Theorem \ref{thm:prob_main} can be turned into an algorithmic procedure with steps given in Algorithm \ref{alg:main}.  The algorithm, which takes as input $\sT$, $q_{0:N}$, $\plant{k\mid  k-1}$ and $c_k(\cdot)$ and outputs  $\left\{\optimalcontrol{k\mid  k-1}\right\}_{1:N}$,  computes $\optimalcontrol{k\mid  k-1}$ following \eqref{eqn:optimal_solution_statement} and  (\ref{eqn:backward_recursion}).  
\begin{algorithm}[h!]
	\caption{Pseudo-code from Theorem \ref{thm:prob_main}: forward control algorithm} \label{alg:main}
	\begin{algorithmic}
			\State \textbf{Inputs:}  $\sT$, $q_{0:N}$, $\plant{k\mid  k-1}$ and $c_k(\cdot)$
		\State \textbf{Output:} $\left\{\optimalcontrol{k\mid  k-1}\right\}_{1:N}$
		\For{$ k = N$  to $1$}
				\State 		Compute $\hat{c}_{k}(\bv{x}_{k})$ using \eqref{eqn:backward_recursion} 
		\State $\bar{c}_{k}(\bv{x}_{k}) \gets c_{k}(\bv{x}_{k}) -\hat{c}_{k}(\bv{x}_{k})$
\State Compute $\bar{p}_{k\mid  k-1}^{(u)}$ using \eqref{eqn:p_bar}
\State $\optimalcontrol{k\mid  k-1} \gets
 \frac{\bar{p}_{k\mid  k-1}^{(u)}\exp\left(-\E_{\plant{k\mid  k-1}}\left[\bar{c}_{k}(\bv{X}_{k})\right]\right)}{\sum_{\bv{u}_{k}}\bar{p}_{k\mid  k-1}^{(u)}\exp\left(-\E_{\plant{k\mid  k-1}}\left[\bar{c}_{k}(\bv{X}_{k})\right]\right)}$
		\EndFor 
	\end{algorithmic}
\end{algorithm}
\begin{Remark}\label{rem:cont}
For continuous variables, the statement of Theorem \ref{thm:prob_main} remains unchanged with the key difference being in the fact that computing the KL-divergence and expectations requires, in general, integrations (and not summations). As a result,  {while} the steps shown in Algorithm \ref{alg:main} remain unchanged when the pfs are continuous{, the integrals need to be estimated and this in general can be computationally difficult.
Next, we discuss a case where the policy can be efficiently computed.}
\end{Remark}
{Consider the setting with continuous variables where  $c_{k}(\bv{x}_{k})=0.5(\bv{x}_{k}-\bv{x}_{d})^{T}\bv{W}(\bv{x}_{k}-\bv{x}_{d})$ and where $\plant{k|k-1}=\mathcal{N}(\bv{A}\bv{x}_{k-1}+\bv{B}\bv{u}_{k},\boldsymbol{\Sigma})$, $\idealplant{k|k-1}=\mathcal{N}(\bv{x}_{d},\bv{R})$,  $\idealcontrol{k|k-1}=\mathcal{N}(\bv{u}_{d},\bv{Q})$, with $\bv{W}\in\R^{n\times n}$, $\boldsymbol{\Sigma}\in\R^{n\times n}$, $\bv{R}\in\R^{n\times n}$ and $\bv{Q}\in\R^{p\times p}$  being positive definite matrices and with $\bv{x}_{d}\in\mathcal{X}$ and $\bv{u}_{d}\in\mathcal{U}$ being constant vectors. An analytical solution when $\bv{W}=0$ was found in \cite{Pegueroles_G+Russo_G_ECC19_confid}. Following similar arguments,  the optimal policy in \eqref{eqn:optimal_solution_statement} - \eqref{eqn:backward_recursion}  becomes
\begin{equation}\label{eqn:gaussian_policy}
\optimalcontrol{k|k-1}=\mathcal{N}(\boldsymbol{\mu}^{\star}_{k},\boldsymbol{\Sigma}_k^\star)
\end{equation}
with $\boldsymbol{\mu}^{\star}_{k}$ and  $\boldsymbol{\Sigma}_k^\star$ given by the following backward recursion for $k=1:N-1$ (starting with $\bv{S}_{N} = 0$)
\begin{equation}\label{eqn:gaussian_policy_recursion}
    \begin{aligned}
    & \bv{S}_{k} = \bv{A}^{T}\left(\bar{\bv{S}}_{k+1}-\bar{\bv{S}}_{k+1}\bv{B}(\bv{Q}^{-1}+\bv{B}^{T}\bar{\bv{S}}_{k+1}\bv{B})^{-1}\bv{B}^{T}\bar{\bv{S}}_{k+1}\right)\bv{A},\\
    & \bar{\bv{S}}_{k}=\bv{S}_{k}+\bv{R}^{-1}+\bv{W},\\
        &\boldsymbol{\Sigma}^\star_{k} = (\bv{Q}^{-1}+\bv{B}^{T}\bar{\bv{S}}_{k}\bv{B})^{-1},\\
        &\boldsymbol{\mu}^{\star}_{k}=-\boldsymbol{\Sigma}_{k}^\star \bv{B}^{T}\bar{\bv{S}}_{k}\bv{A}\bv{x}_{k-1}+\boldsymbol{\Sigma}^\star_{k}(\bv{B}^{T}\bv{\bar{S}}_{k}\bv{x}_{d}+\bv{Q}^{-1}\bv{u}_{d}).
    \end{aligned}
\end{equation}
See the appendix for the derivations.}
\paragraph*{{\bf Running Example (continue).}}
We now use Theorem \ref{thm:prob_main} (and hence Algorithm \ref{alg:main}) to swing-up the target pendulum introduced in the previous part of the example.  The cost associated to this control task is given in 
\eqref{eqn:state_cost}. For reasons that will be clear later,  for both the continuous and discrete settings described above,  at each $k$ we sampled the control input from the policy computed (given $\bv{x}_{k-1}$) via Algorithm \ref{alg:main} with $N=1$.  As shown in Figure \ref{fig:sim_forward_pendulum}, in both the continuous and discrete settings, the policy $\optimalcontrol{k \mid k-1}$ computed via Algorithm \ref{alg:main} stabilized the target pendulum. In the figure, the behavior is shown for the angular position of the controlled pendulum when the control input is sampled, at each $k$, from $\optimalcontrol{k \mid k-1}$. See our github for the implementation details at \url{https://tinyurl.com/46uccxsf}.
~

\subsection{Tackling the Inverse Control Problem}\label{sec:learning_cost}

First, we give a result to reconstruct, via a convex problem, $\bar{c}_k(\cdot)$ by observing a sequence of states sampled from $\plant{k\mid  k-1}$ when the control inputs sampled from $\optimalcontrol{k\mid k-1}$ are applied. Then, we build on  this result to reconstruct  ${c}_k(\cdot)$.  {We let Assumption \ref{asn:inverse_problem} hold and let  $\{(\hat{\bv{x}}_{k-1},\hat{\bv{u}}_{k})\}_{1:M}$ be the observed data pairs. We consider the likelihood function of $\bv{w}= [\bv{w}_1^T,\ldots \bv{w}_M^T]^T$ 
\begin{align}\label{eqn:log-likelihood}
\begin{split}
& L\left({\bv{w}}\right):=  \prod_{k=1}^M  \\
& \frac{\modifiedcontrol{k\mid  k-1}(\hat{\bv{x}}_{k-1},\hat{\bv{u}}_k) \exp\left(
\sum_{\bv{x}_k}p_{{k}}(\bv{x}_k\mid\hat{\bv{x}}_{k-1},\hat{\bv{u}}_k)\bv{w}_k^T{\bv{h}({{\bv{x}}}_k})\right)}{\sum_{\bv{u}_{k}}\modifiedcontrol{k\mid  k-1}(\hat{\bv{x}}_{k-1},{\bv{u}}_k) \exp\left(
\sum_{\bv{x}_k}p_{{k}}(\bv{x}_k\mid\hat{\bv{x}}_{k-1},{\bv{u}}_k)\bv{w}_k^T{\bv{h}({{\bv{x}}}_k})\right)},\\
\end{split}
\end{align}
where 
\begin{equation}\label{eqn:modified_control}
\begin{split}
    &\modifiedcontrol{k\mid  k-1}({\bv{x}}_{k-1},{\bv{u}}_k):=q({\bv{u}}_k\mid{\bv{x}}_{k-1})\\ 
    &\exp\left(-\DKL\left(p(\bv{x}_k\mid{\bv{x}}_{k-1},{\bv{u}}_k)\mid\mid q(\bv{x}_k\mid{\bv{x}}_{k-1},{\bv{u}}_k)\right)\right).
\end{split}
\end{equation}}
{and where we used the shorthand notation} $p_{{k}}(\bv{x}_k\mid\hat{\bv{x}}_{k-1},{\hat{\bv{u}}}_k)$,  $p_{{k}}(\bv{x}_k\mid\hat{\bv{x}}_{k-1},{{\bv{u}}}_k)$ and $q_{{k}}({\bv{u}}_k\mid\hat{\bv{x}}_{k-1})$ to denote $p^{(x)}_{{k}}(\bv{x}_k\mid\bv{x}_{k-1} = \hat{\bv{x}}_{k-1},\bv{u}_k = {\hat{\bv{u}}}_k)$,  $p^{(x)}_{{k}}(\bv{x}_k\mid\bv{x}_{k-1} = \hat{\bv{x}}_{k-1},\bv{u}_k)$ and $q^{(u)}_{{k}}({\bv{u}}_k\mid\bv{x}_{k-1}= \hat{\bv{x}}_{k-1})$, respectively.  {The next result is related to what was independently reported in \cite{IRL-Todorov} in a different context (see Remark \ref{rem:Todorov}) and gives the maximum likelihood estimate for $\bv{w}$.}

\begin{Theorem}\label{thm:estimator}
Let Assumption \ref{asn:inverse_problem} hold and let the observed data {pairs} be $\{(\hat{\bv{x}}_{k-1},\hat{\bv{u}}_{k})\}_{1:M}$, 
with $\hat{\bv{x}}_{k} \sim \plant{k\mid  k-1}$, $\hat{\bv{u}}_{k} \sim \optimalcontrol{k\mid k-1}$ and where $\optimalcontrol{k\mid k-1}$ is  from Theorem \ref{thm:prob_main}. 
 Then, {the maximum likelihood estimate for $\bv{w}$ with likelihood function \eqref{eqn:log-likelihood} is} obtained by solving the convex  problem
\begin{equation}\label{eqn:prob_log_likelihood_statement}
    \begin{aligned}
    &{\bv{w}}^\star:=\left[\bv{w}^{\star T}_1,\ldots,\bv{w}^{\star T}_M\right]^T \in \\
    \underset{{\bv{w}}}{\text{arg min}}
     &\Bigg\{\sum_{k=1}^M\bigg({-\E_{p(\bv{x}_k\mid\hat{\bv{x}}_{k-1},{\hat{\bv{u}}}_k)}\left[{\bv{w}}_k^T{\bv{h}({\bv{x}}}_k)\right]} \\
     & +\ln\Big(\sum_{\bv{u}_k}\modifiedcontrol{k\mid  k-1}(\hat{\bv{x}}_{k-1},{\bv{u}}_k) \\ &\exp\left(\E_{p(\bv{x}_k\mid\hat{\bv{x}}_{k-1},\bv{u}_k)}\left[{\bv{w}}_k^T{\bv{h}({{\bv{x}}_k}})\right]\right)\Big)\bigg)\Bigg\},
    \end{aligned}
\end{equation}
where {$\modifiedcontrol{k\mid  k-1}(\cdot,\cdot)$ is defined in \eqref{eqn:modified_control}}.
\end{Theorem}
{Before giving the proof, we make the following
\begin{Remark} In Theorem \ref{thm:estimator},  $M$,  the number of observed data pairs $(\hat{\bv{x}}_{k-1},\hat{\bv{u}}_{k})$,  can be different from the time horizon, $N$, of the forward control problem (Problem \ref{prob:main}).  This is useful when, as in e.g.,  our running example,  observed data are obtained from multiple experiments.
\end{Remark}}
\begin{pf}
Assumption \ref{asn:inverse_problem}, together with the fact that the observed actions are sampled from the policy of Algorithm \ref{alg:main} implies that sequence of control inputs $\{\hat{\bv{u}}_k\}_{1:M}$ are determined by sampling, at each $k$, from the  pmf
\begin{equation}
\begin{aligned}
&\optimalcontrol{k\mid k-1}= \\& \frac{\modifiedcontrol{k\mid k-1}(\hat{\bv{x}}_{k-1}, {\bv{u}}_k)\exp\left(\E_{p_k(\bv{x}_k \mid \hat{\bv{x}}_{k-1}, \bv{u}_k)} \left[\bv{w}_k^T {\bv{h}({\bv{x}}_k})\right]\right)}
{\sum_{\bv{u}_{k}} \modifiedcontrol{k\mid k-1}(\hat{\bv{x}}_{k-1}, {\bv{u}}_k)\exp\left(\E_{p_k(\bv{x}_k \mid \hat{\bv{x}}_{k-1}, \bv{u}_k)}\left[ \bv{w}_k^T {\bv{h}({\bv{x}}_k})\right]\right)}{.}
\end{aligned}
\end{equation}
{This leads to the likelihood function in \eqref{eqn:log-likelihood}.  We find 
$$
\bv{w}^\star \in \underset{\bv{w}}{\text{arg} \max} L(\bv{w})
$$
by minimizing the negative log-likelihood $-\log L(\bv{w})$. In turn, } the negative log-likelihood can be written as 
\begin{equation}\label{eqn:log-likelihood_main}
\begin{split}
& \sum_{k=1}^M({-\ln(\modifiedcontrol{k\mid  k-1}(\hat{\bv{x}}_{k-1},\hat{\bv{u}}_k))}\\
&-\sum_{\bv{x}_k}p_{{k}}(\bv{x}_k\mid\hat{\bv{x}}_{k-1},\hat{\bv{u}}_k)\bv{w}_k^T{\bv{h}({{\bv{x}}}_k})) \\
& + \sum_{k=1}^M\ln(\sum_{\bv{u}_k}\modifiedcontrol{k\mid  k-1}(\hat{\bv{x}}_{k-1},\bv{u}_k)\cdot\\
& \exp(\sum_{\bv{x}_k}p_{{k}}(\bv{x}_k\mid\hat{\bv{x}}_{k-1},{\bv{u}}_k)\bv{w}_k^T{\bv{h}({{\bv{x}}}_k})))
\end{split}
\end{equation}
In \eqref{eqn:log-likelihood_main} the term $\modifiedcontrol{k\mid  k-1}(\cdot,\cdot)$ is defined in (\ref{eqn:modified_control}). Hence, minimizing the negative log-likelihood results into the following unconstrained optimization problem:
 \begin{equation}\label{eqn:problem_log_likelihood_first}
    \begin{aligned}
    \underset{{\bv{w}}}{\text{min}}& {-
    \log}L(\bv{w}),
    \end{aligned} 
\end{equation}
where ${-
    \log}L(\bv{w})$ is given in \eqref{eqn:log-likelihood_main}. Moreover, the first term in the cost function does not depend on the decision variable $\bv{w}${.} Hence, such term can be removed from the formulation of the problem and this yields (\ref{eqn:prob_log_likelihood_statement}).  Problem's convexity follows from the cost function being a convex combination of a linear function and the log-sum-exponential function (which is convex). \hspace*{\fill}\qed
\end{pf}
{\begin{Remark}\label{rem:Todorov}
Following Theorem \ref{thm:estimator}, the estimated $\bar{c}_k\left({\bv{x}}_k\right)$, say $\bar{c}_k^\star\left({\bv{x}}_k\right)$, is given by $\bar{c}_k^\star\left({\bv{x}}_k\right) =   -{{\bv{w}}^\star_k}^T{\bv{h}(\bv{x}}_k)$.  As in \cite{IRL-Todorov} in Theorem \ref{thm:estimator} we exploit the fact that the optimal solution from Theorem \ref{thm:prob_main} belongs to the exponential family.  In \cite{IRL-Todorov} the framework of linearly solvable MDPs is exploited, which relies on the assumption that $p^{(x)}_k	\left(	\mathbf{x}_k \mid \mathbf{u}_k, \mathbf{x}_{k-1} \right)p^{(u)}_k	\left(	\mathbf{u}_k \mid  \mathbf{x}_{k-1} \right) = p	\left(	\mathbf{x}_k \mid \mathbf{x}_{k-1} \right)$ and $q^{(x)}_k	\left(	\mathbf{x}_k \mid \mathbf{u}_k, \mathbf{x}_{k-1} \right)q^{(u)}_k	\left(	\mathbf{u}_k \mid  \mathbf{x}_{k-1} \right) = q	\left(	\mathbf{x}_k \mid \mathbf{x}_{k-1} \right)$.  This assumption is not required in our result.
\end{Remark}}
\begin{Remark}\label{cons_wi}
The problem in (\ref{eqn:prob_log_likelihood_statement}) is an unconstrained convex optimization problem with a twice differentiable cost.  Constraints on the $\bv{w}_k$'s can be added to the problem to capture additional desired properties of the cost. 
\end{Remark}

Consider the situation where in Theorem \ref{thm:estimator} at each $k$ the policy $\optimalcontrol{k\mid k-1}$ is  obtained from Theorem \ref{thm:prob_main} with $N=1$. In this case,  the result provides an estimate of the cost $c_k(\cdot)$.  Next, we consider the case where the cost, which we simply denote by $c(\cdot)$, is stationary. 

\begin{Corollary}\label{clry: estimator}
Let Assumption \ref{asn:inverse_problem} hold,  the cost be stationary and $\optimalcontrol{k\mid k-1}$ be the policy obtained at each $k$ from Theorem \ref{thm:prob_main} with $N=1$.  Then,  ${c}^\star(\bv{x}_{k}) = -{\bv{w}_{s}^\star}^T{\bv{h}(\bv{x}}_k)$,  where $\bv{w}_{s}^\star$ is given by:
  ~\begin{equation}\label{eqn:estimation_stationary}
{\begin{aligned}
\bv{w}_{s}^{\ast} \in \underset{\bv{w}_{s}}{\text{arg min}} & \Bigg\{\sum_{k=1}^M\left(
-\E_{p(\bv{x}_k\mid\hat{\bv{x}}_{k-1},{\hat{\bv{u}}}_k)}\left[\bv{w}_{s}^
T{\bv{h}({{\bv{x}}}_k})\right]\right )\\
 & +\sum_{k=1}^M\ln\left(\sum_{\bv{u}_k}\modifiedcontrol{k\mid  k-1}(\hat{\bv{x}}_{k-1},\bv{u}_k)\right.\\
&\exp\left(\E_{p(\bv{x}_k\mid\hat{\bv{x}}_{k-1},{\bv{u}}_k)}\left[\bv{w}_{s}^T{\bv{h}({{\bv{x}}}_k})\right]\right)\bigg)\Bigg\}.
\end{aligned}}
\end{equation}
 with
$\bv{w}_{s} \in \mathbf{R}^{F}$ and $\modifiedcontrol{k\mid  k-1}(\hat{\bv{x}}_{k-1},\bv{u}_k)$ defined as in Theorem \ref{thm:estimator}.
\end{Corollary}

\begin{pf}
Follows from Theorem \ref{thm:estimator} after noticing that, when $N=1$ and the cost is stationary, Theorem \ref{thm:prob_main}  yields
$$\optimalcontrol{k\mid  k-1} =
 \frac{\bar{p}_{k\mid  k-1}^{(u)}\exp\left(-\E_{\plant{k\mid  k-1}}\left[{c}(\bv{X}_{k})\right]\right)}{\sum_{\bv{u}_{k}}\bar{p}_{k\mid  k-1}^{(u)}\exp\left(-\E_{\plant{k\mid  k-1}}\left[{c}(\bv{X}_{k})\right]\right)}.  \hspace*{\fill}\qed
 $$
\end{pf}

Following Corollary \ref{clry: estimator},  when the cost is stationary, one needs to solve an optimization problem that has as decision variable the $\bv{w}_{s}\in\mathbf{R}^F$ rather than $\bv{w}\in\mathbf{R}^{ MF}$. Also, Corollary \ref{clry: estimator} only requires that data are collected via a {\em greedy} policy obtained at each $k$ from Theorem \ref{thm:prob_main} with $N=1$.  {T}his is convenient as it bypasses the need to solve the backward recursion in \eqref{eqn:backward_recursion}.

\subsubsection{Turning Corollary \ref{clry: estimator} into an Algorithm}
We turn Corollary \ref{clry: estimator} into an algorithmic procedure with steps given in Algorithm \ref{alg:estimator}.  An algorithm for Theorem \ref{thm:estimator} can also be obtained and it is omitted here for brevity.

\begin{algorithm}[ht]
	\caption{Pseudo-code from Corollary \ref{clry: estimator}: inverse control algorithm} \label{alg:estimator}
	\begin{algorithmic}
			\State \textbf{Inputs:}  observed data $\{\hat{\bv{u}}_k\}_{1:M}$ and $\{\hat{\bv{x}}_{k-1}\}_{1:M}$, 
            \State \hspace{1.1cm} $F$-dimensional features vector $\bv{h}(\bv{x}_k)$, 
            \State \hspace{1.1cm} $\plant{k \mid k-1}$,  $\idealplant{k \mid k-1}$, $\idealcontrol{k \mid k-1}$ 
		\State \textbf{Output:} ${c}^\star\left({\bv{x}}_k\right)$
		\For{$ k = 1$  to $M$}
            \State Compute $\modifiedcontrol{k\mid  k-1}(\hat{\bv{x}}_{k-1},{\bv{u}}_k)$ using \eqref{eqn:modified_control}
            \EndFor
        \State Compute $\bv{w}_{s}^\star$ by solving the problem in \eqref{eqn:estimation_stationary}
        \State ${c}^\star(\bv{x}_{k}) \gets -{\bv{w}_{s}^\star}^T{\bv{h}(\bv{x}}_k)$
	\end{algorithmic}
\end{algorithm}

\subsection{Special Case}\label{sec:special_cases}

We now discuss a special case of our results, relevant for applications, when the reference pfs are uniform.

\noindent{\bf Forward problem.} The KL-divergence component in the cost of Problem \ref{prob:main} becomes an entropic regularizer  and the optimal policy in \eqref{eqn:optimal_solution_statement} then becomes
\begin{equation}\label{eqn:optimal_solution_statement_uniform}
\begin{split}
& \optimalcontrol{k\mid  k-1} =\\
&\frac{\exp\left(-\E_{\plant{k\mid  k-1}}\left[\ln\left(\plant{k|k-1}\right)+\bar{c}_{k}(\bv{X}_{k})\right]\right)}{\sum_{\bv{u}_{k}}\exp\left(-\E_{\plant{k\mid  k-1}}\left[\ln\left(\plant{k|k-1}\right)+\bar{c}_{k}(\bv{X}_{k})\right]\right)},
\end{split}
\end{equation} 
where $\bar{c}_{k}(\bv{x}_{k})= c_{k}(\bv{x}_{k}) -\hat{c}_{k}(\bv{x}_{k})$ and with the backward recursion in \eqref{eqn:backward_recursion} becoming
\begin{equation*}
\begin{split}
&\hat{c}_{k}(\bv{x}_{k})= \ln\Bigg(\sum_{\bv{u}_{k}}\exp\bigg(-\E_{\plant{k+1\mid  k}}\Big[\ln\left(\plant{k+1|k}\right)\\
& +\bar{c}_{k+1}(\bv{X}_{k+1})\Big]\bigg)\Bigg),\\
&\E_{\plant{N+1\mid  N}}\left[\ln\left(\plant{N+1|N}\right) + \bar{c}_{N+1}(\bv{X}_{N+1})\right] = 0.
\end{split}
\end{equation*}

\noindent{\bf Inverse problem.} Since in the optimal policy is now the one in \eqref{eqn:optimal_solution_statement_uniform},  the problem  in \eqref{eqn:prob_log_likelihood_statement} becomes
\begin{equation*}
    \begin{aligned}
    &\underset{{\bv{w}}}{\text{arg min}}
     \Bigg\{\sum_{k=1}^M\Bigg({-\E_{p(\bv{x}_k\mid\hat{\bv{x}}_{k-1},{\hat{\bv{u}}}_k)}\left[{\bv{w}}_k^T{\bv{h}({\bv{x}}}_k)\right]} \\
     &+\ln\bigg(\sum_{\bv{u}_k}\exp\Big(\E_{p(\bv{x}_k\mid\hat{\bv{x}}_{k-1},{\bv{u}}_k)}\Big[-\ln (p(\bv{x}_k\mid\hat{\bv{x}}_{k-1},{\bv{u}}_k))\\
     &+  {\bv{w}}_k^T{\bv{h}({{\bv{x}}_k}})\Big]\Big)\bigg)\Bigg)\Bigg\}.
    \end{aligned}
\end{equation*}

\begin{figure*}
    \centering
    \includegraphics[width=0.25\linewidth]{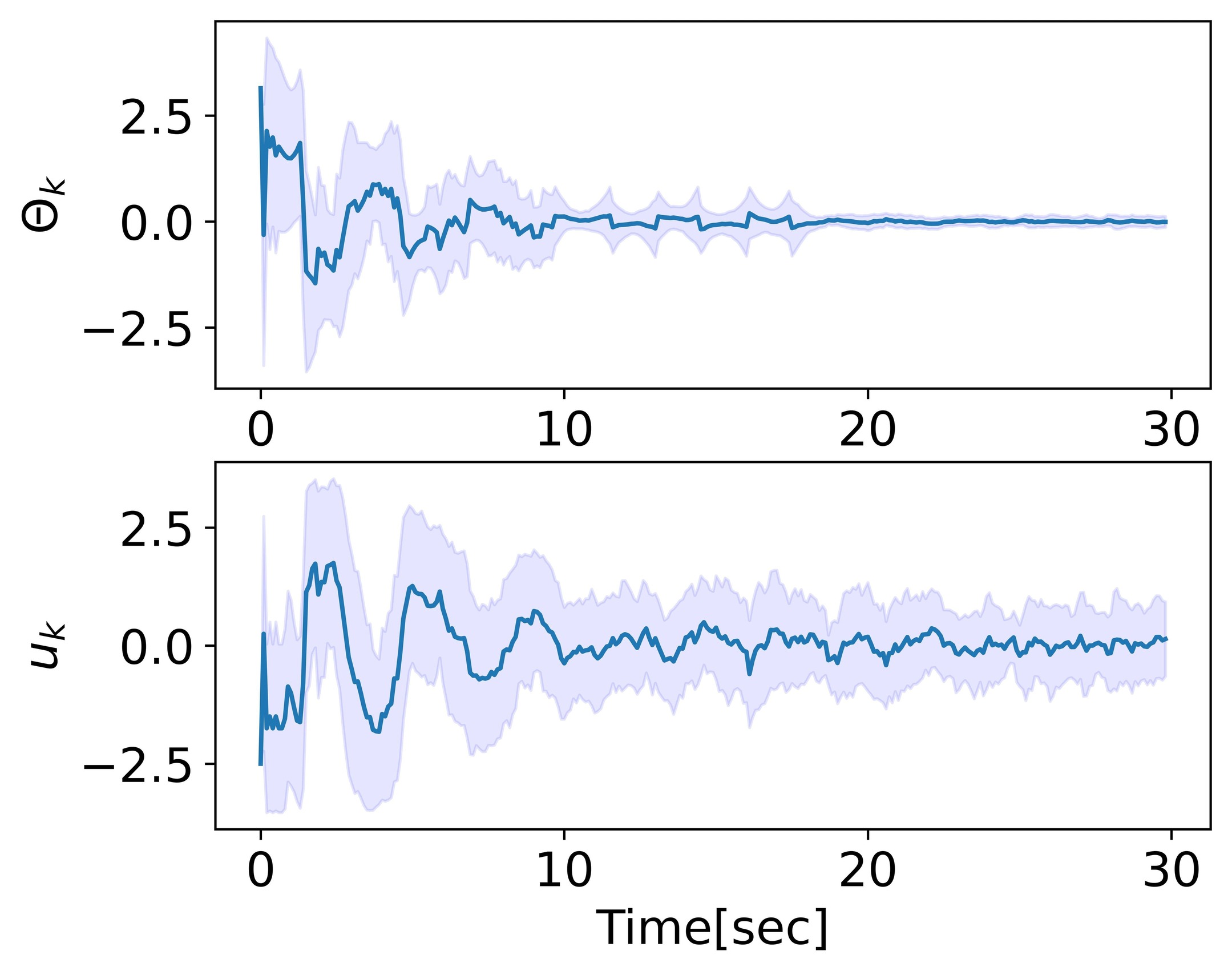}
        \includegraphics[width=0.25\linewidth]{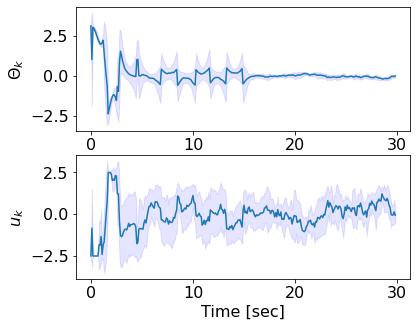}
    \includegraphics[width=0.25\linewidth]{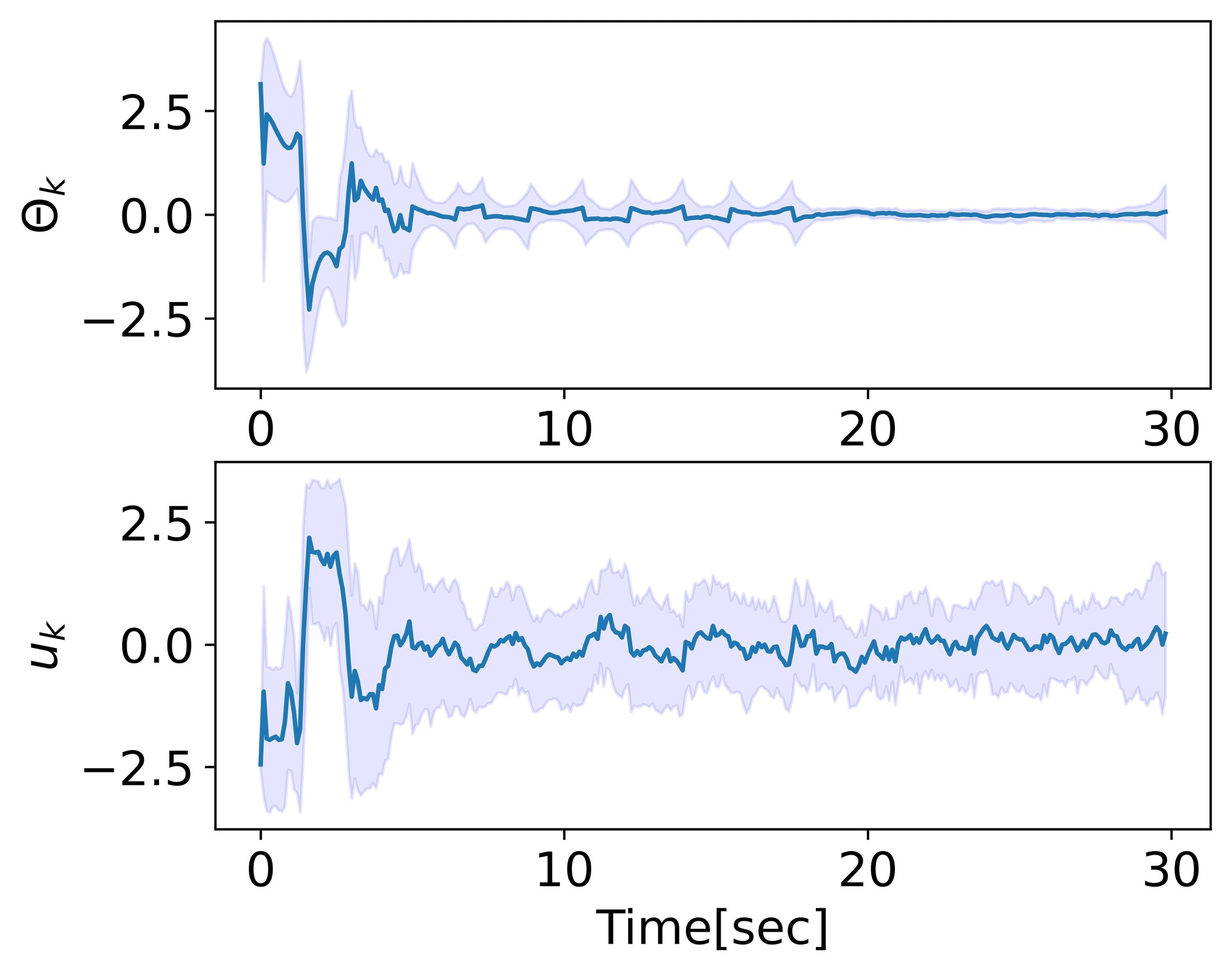}
    \caption{Angular position and  control input of the target Pendulum when the pf is estimated via the histogram filter (left{and middle} panels) and Gaussian Processes (right panels).  Figures  obtained from $20$ simulations, using ${c}^{\star}(\cdot)$ as an input to Algorithm \ref{alg:main}. Bold lines represents the mean;  the shaded region is confidence interval corresponding to the standard deviation.}
    \label{fig:inverse_control_pendulum}
\end{figure*}

\begin{figure}[b!]
    \centering
    \includegraphics[width=0.4\columnwidth]{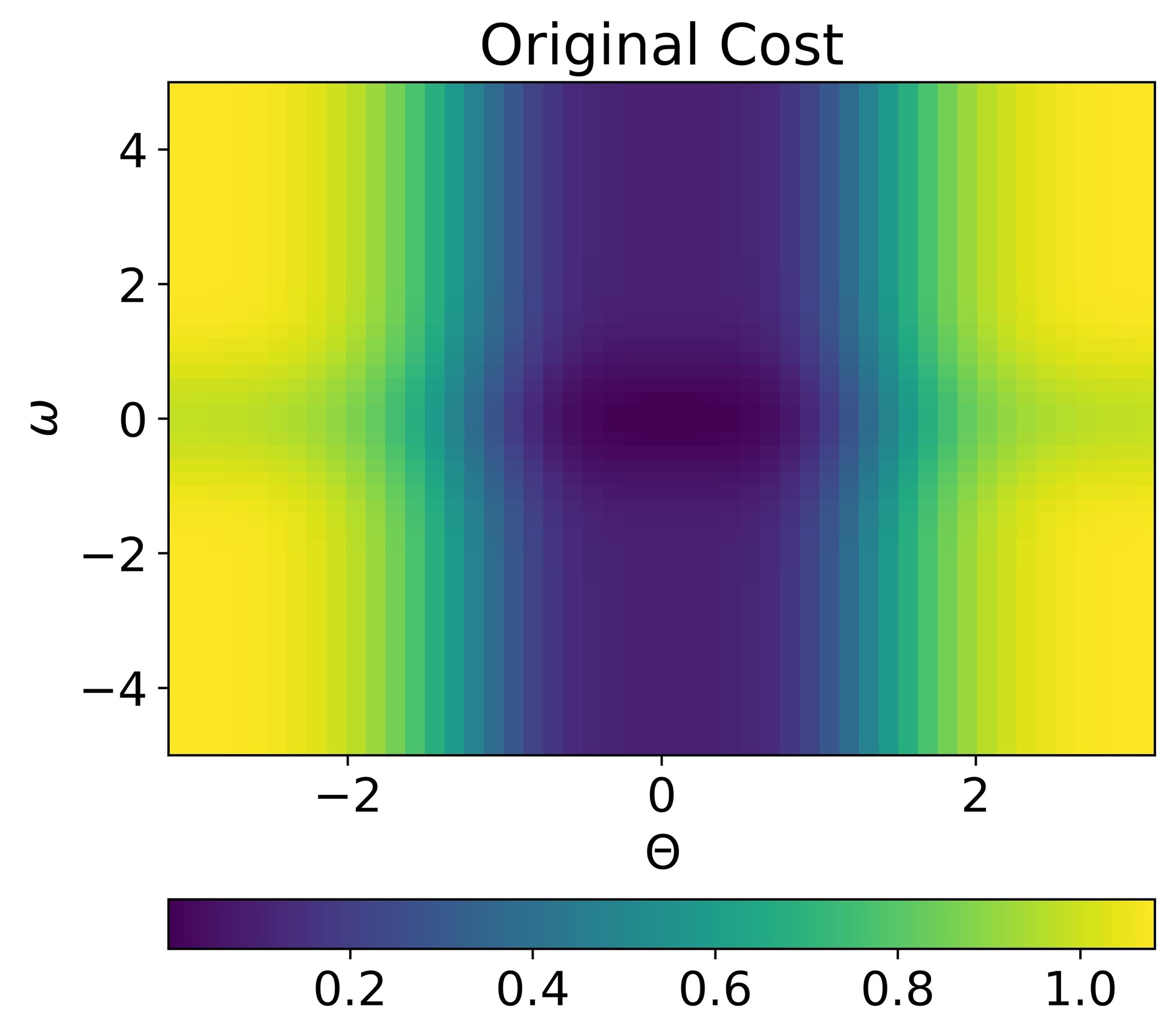}
    \includegraphics[width=0.4\columnwidth]{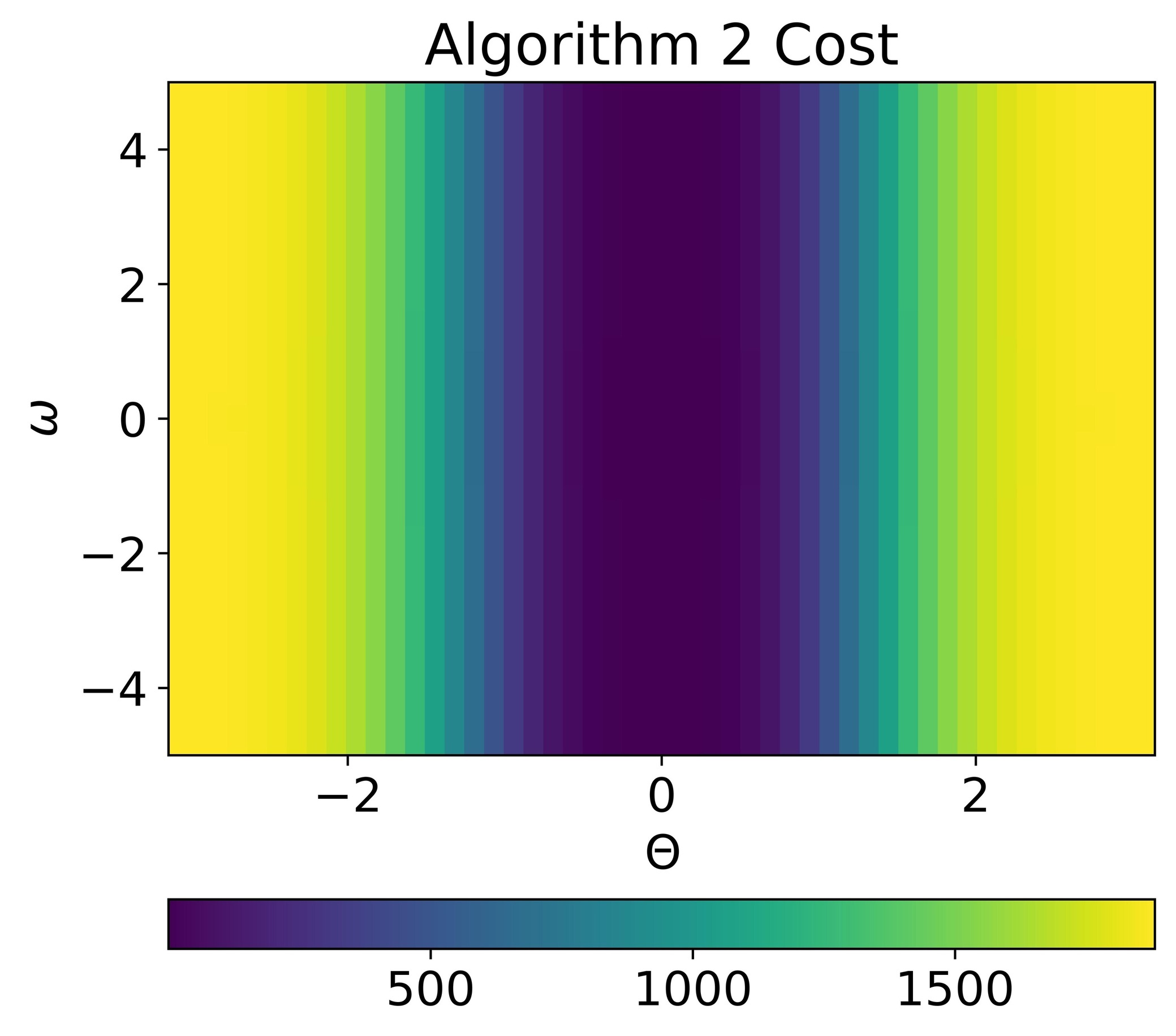}
    \includegraphics[width=0.4\columnwidth]{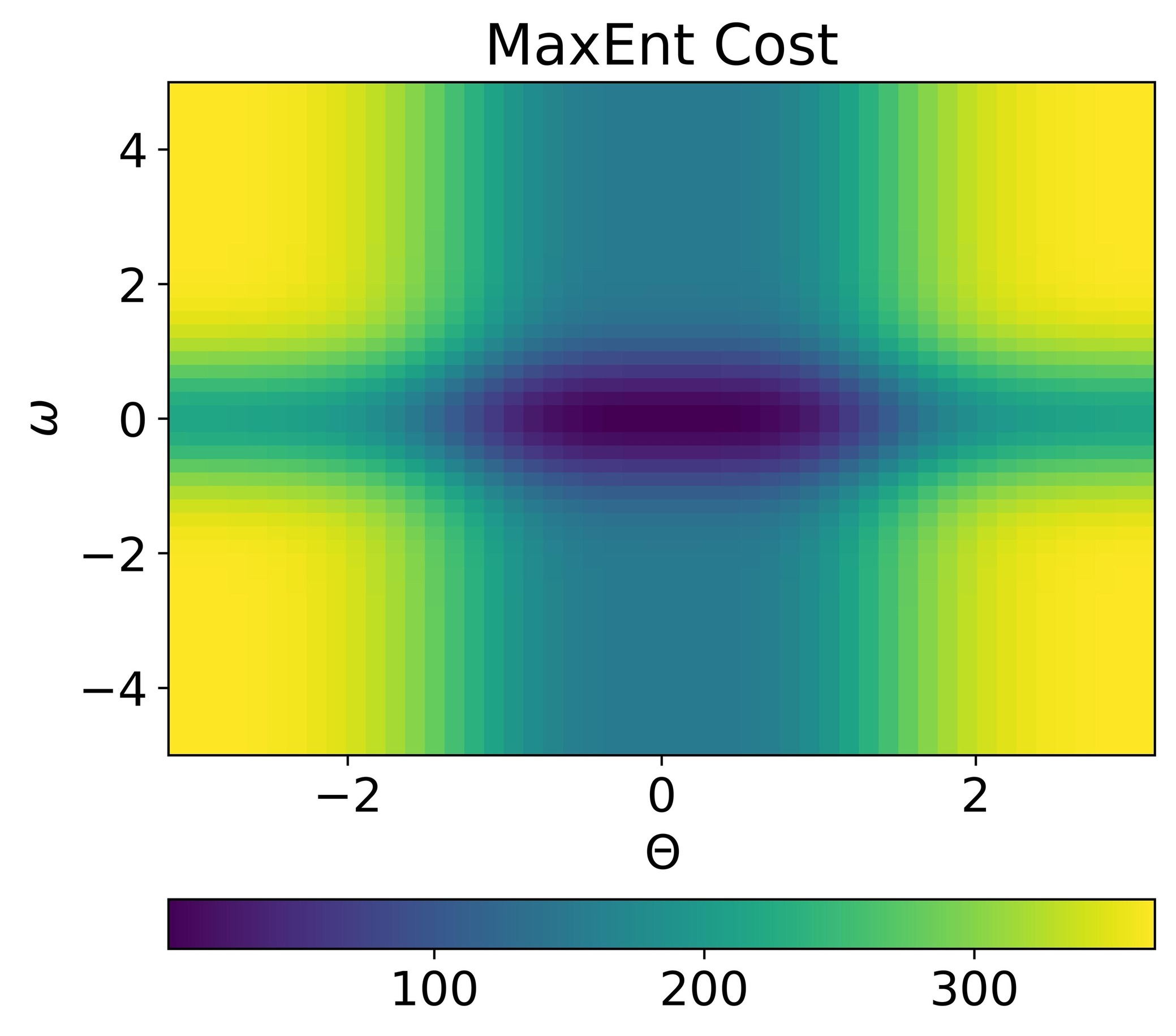}
    \includegraphics[width=0.4\columnwidth]{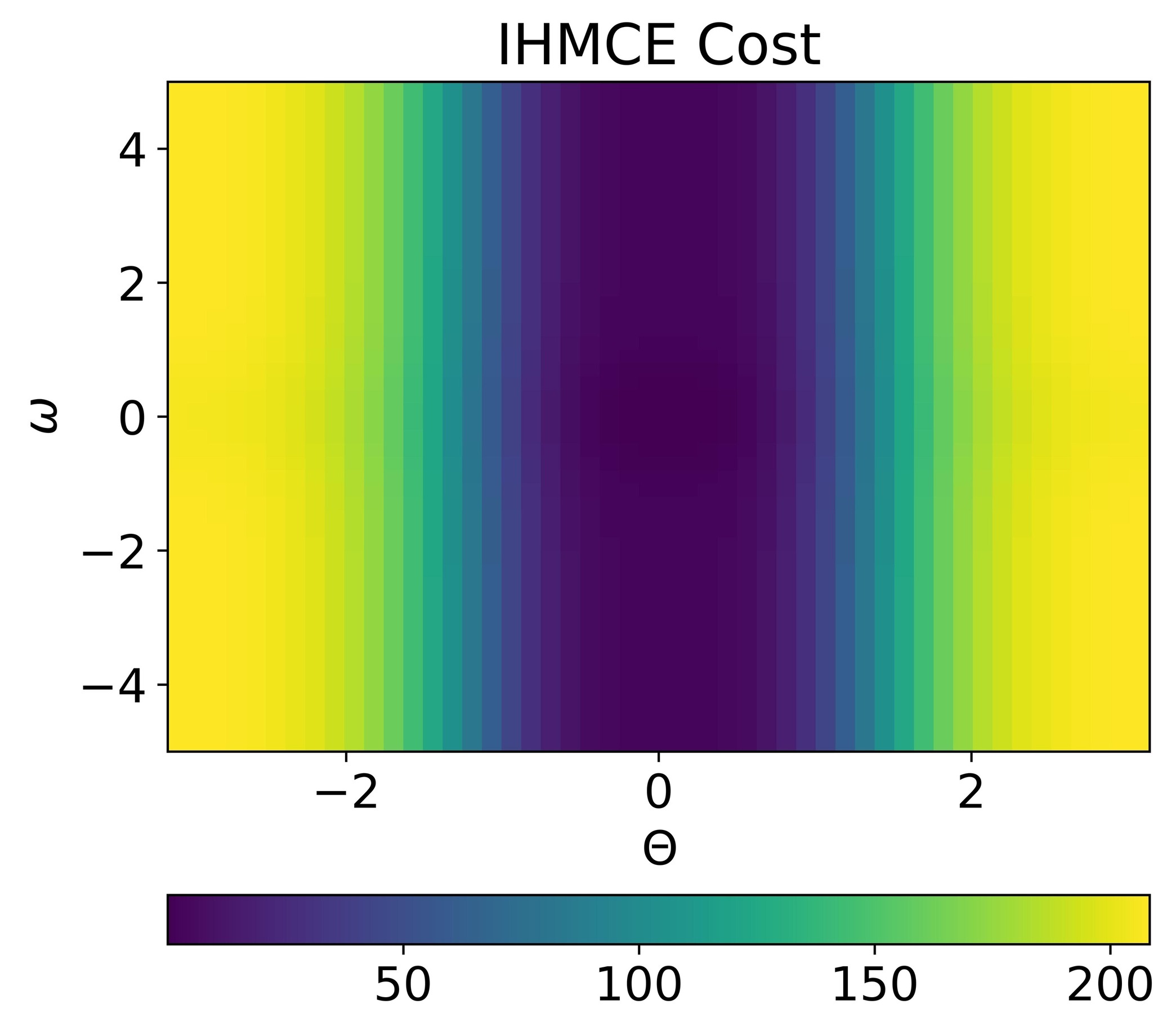}
    \caption{{Top left: original cost function. In the other panels the cost reconstructed via: Algorithm \ref{alg:estimator} (top-right),  MaxEnt (bottom-left) and IHMCE (bottom-right).}}
    \label{fig:expert_value}
\end{figure}

\paragraph*{{\bf Running Example (continue).}}
In this final part of the example we show the application of Algorithm \ref{alg:estimator} when this is leveraged to reconstruct the cost  in \eqref{eqn:state_cost} that was given as an input to Algorithm \ref{alg:main}.   Continuous and discrete settings are discussed separately.  {For the discrete setting we  give a  comparison with \cite{10.5555/1620270.1620297,8115277}.} See \url{https://tinyurl.com/46uccxsf} for  details. \\
\\
\noindent {\bf Discrete setting.} We used approximately $300$ data {pairs} collected from a single simulation of the target pendulum controlled via the policy computed in the previous part of the running example (from this dataset,  we removed the points where there were discontinuities due to angle wrapping).  These were the observed data given as an input to Algorithm \ref{alg:estimator}.  Also, we defined the features as $\bv{h}(\bv{x}_{k})= \left[\vert\theta_{k}-\theta_{d}\vert,\vert\omega_{k}-\omega_{d}\vert\right]^T$ and obtained from Algorithm \ref{alg:estimator} (using CVX \cite{diamond2016cvxpy} to solve the optimization problem) the weights $\bv{w}_{s}^{\star} =[-6.09, -4.48]^{T}$.  Hence, in accordance with Corollary \ref{clry: estimator},  we obtained ${c}^{\star}(\bv{x}_{k}) = 6.09 \vert\theta_{k}-\theta_{d}\vert + 4.48 \vert\omega_{k}-\omega_{d}\vert$.  Then,  we used this cost as input for Algorithm \ref{alg:main} (again with $N=1$) to solve the forward control problem. As shown in Figure \ref{fig:inverse_control_pendulum} (left panel) the policy computed using ${c}^{\star}(\cdot)$  effectively stabilizes the pendulum.   {To further validate our results, we used Algorithm \ref{alg:estimator} with the same dataset described above but with $\bv{h}(\bv{x}_{k})=[1-\exp(-(\cos(\theta_{k})-1)^{2}),1-\exp(-\omega_{k}^{2})]$.  This time we performed $10$ experiments. For each experiment, the algorithm provided a set of weights and the time taken by Algorithm \ref{alg:estimator} to output the weights is summarized in Table \ref{tab:comparison} (first line).  We also quantified the similarity between the original cost and the reconstructed cost.  To this aim, we normalized both the original cost and the cost reconstructed via Algorithm \ref{alg:estimator} and then used the KL divergence to quantify their discrepancy.  This measure of discrepancy, which is averaged across the $10$ experiments, is also shown in Table \ref{tab:comparison} (last column).  In Figure \ref{fig:expert_value} (top-left) the reconstructed cost is shown for a set of randomly selected weights across the experiments ($\bv{w}_{s}^{\star} =[-1913.3,-3.2]^{T}$).  Finally,  given this set of weights we used the corresponding reconstructed cost as input for Algorithm \ref{alg:main} and the results are shown in Figure \ref{fig:inverse_control_pendulum} (middle panel).}
\begin{table}[t!]
\centering
\caption{{Time taken by Algorithm \ref{alg:estimator} to reconstuct the cost and comparison with  Max-Ent (VI) and IHMCE (Soft-VI and MCS).  The table also reports the discrepancy between the original and the reconstructed costs.  Results obtained from $10$ experiments.  Bold: the best values for each column.}}
\resizebox{\columnwidth}{!}{%
{\begin{tabular}{p{4cm}cccc}
\toprule
Method & Best Case & Worst Case & Average &  Discrepancy \\
\midrule
Algorithm \ref{alg:estimator} & {\bf 32.2 sec} & {\bf 1.05 min} & {\bf 40 sec} & {\bf 0.00162} \\
Max-Ent with VI \citep{8115277} & 15 mins & 40 mins & 28 mins & 0.124 \\
IHMCE with soft VI and MCS  \citep{10.5555/1620270.1620297} & 10 mins & 36 mins & 25 mins & 0.00467 \\
\bottomrule
\end{tabular}}
}
\label{tab:comparison}
\end{table}

{\noindent {\bf Benchmarking Algorithm \ref{alg:estimator}.} By using the setting introduced so far, we now present a comparison between the results obtained via Algorithm \ref{alg:estimator} with the Maximum Entropy-IRL algorithm (MaxEnt for short in what follows) from \cite{10.5555/1620270.1620297} and the Infinite Time Horizon Maximum Causal Entropy IRL (IHMCE) from \cite{8115277}.  Both these methods rely on discretized action/state spaces and assume (as also required by Algorithm \ref{alg:estimator}) a linear parametrization in a set of features. Additionally, these methods also assume that the {policy} for which the cost is reconstructed is a soft-max (see Remark \ref{rem:exponential_policy} for a discussion on this point and at the end of Section \ref{sec:contributions}).  The results for MaxEnt and IHMCE were obtained by using the last set of features introduced above and the same data that were used as input to Algorithm \ref{alg:estimator}.  Fully documented code implementing MaxEnt and IHMCE  is also available at \url{https://tinyurl.com/46uccxsf},  together with a test case validating our implementation on the  grid-world from \cite[Chapter $3$]{Sutton:1998:IRL:551283}.  With our first set of experiments we aimed at reconstructing the pendulum cost using MaxEnt, which foresees both forward and backward pass (FP and BP) steps.  Our experiments revealed that, for the pendulum of this running example,  $3.75e^{10}$ iterations were required for a single BP step to converge and this made the approach {infeasible} for the application.  The high computation cost of the BP step was already noted in \cite{8115277} where it was also proposed to use Value Iteration (VI) rather than BP to improve the computational aspect. Hence, inspired by this,  we then implemented a version of MaxEnt with VI.  The average discrepancy between the cost reconstructed with this algorithm and the original cost across $10$ experiments is shown in Table \ref{fig:expert_value}, together with the corresponding computation time (one representative cost, obtained by randomly selecting one of the weight vectors from the $10$ experiments is shown in the bottom-left panel of Figure \ref{fig:expert_value}). Based on these findings, we then implemented the IHMCE algorithm from \cite{8115277} that makes use of soft VI and Monte Carlo Simulations (MCS).  Again, we computed the average discrepancy across $10$ experiments and recorded the computation times. The results are also given in Table \ref{fig:expert_value} and a representative cost is  shown in Figure \ref{fig:expert_value} (bottom-right). In summary,  our numerical studies illustrate that Algorithm \ref{alg:estimator} is both significantly faster than Max-Ent (with VI) and IHMCE (with soft VI and MCS) and, at the same time, achieves lower cost discrepancy.}\\
\\
\noindent {\bf Continuous setting.} Again, we used a dataset of $300$ data-points collected from a single simulation where the target pendulum was controlled by the policy from Algorithm \ref{alg:main} as described in the previous part of the  example.  This time we used as features vector $\bv{h}(\bv{x}_{k})= \left[(\cos(\theta_{k})-\cos(0))^{2},(\cos(\theta_{k})-\cos(\pi))^{2}\right]^T$. When this vector was given as an input to  Algorithm \ref{alg:estimator} (using again CVX to solve the problem) the weights we obtained were $\bv{w}_{s}^{\star} =[-13.66,8.50]^{T}$ so that ${c}^{\star}(\bv{x}_{k}) = 13.66 (\cos(\theta_{k})-\cos(0))^{2} - 8.50 (\cos(\theta_{k})-\cos(\pi))^{2}$. Note that the first term in ${c}^{\star}(\cdot)$ conveniently drives the pendulum towards the unstable equilibrium while the second term is pushing it away from the stable equilibrium in $\pi$, consistently with the cost in \eqref{eqn:state_cost}. Again, we then gave ${c}^{\star}(\cdot)$ as an input to Algorithm \ref{alg:main} and verified that it was in fact able to swing-up the pendulum. The numerical results are reported in Figure \ref{fig:inverse_control_pendulum} (right panel). 
~

\section{Application Example}\label{sec:application}
We use our results to tackle an application involving routing unicycle robots in an environment with obstacles.  {Specifically, we consider two scenarios: in Scenario $1$, the robot moves in an environment without obstacles and needs to achieve a desired destiation; (ii) in Scenario $2$, the robot needs to reach the} destination while avoiding obstacles.  {For both scenarios,} we first use Algorithm \ref{alg:main} to compute a policy routing the robot.  Then, we leverage Algorithm \ref{alg:estimator} to reconstruct the navigation cost from observed robot trajectories.  The results are validated both via simulations and via real hardware experiments,  leveraging the {\em Robotarium} platform that offers both a hardware infrastructure and a high-fidelity simulator  \cite{robotarium2020}.  The robots in the Robotarium are unicycles of dimensions $11 \text{cm} \times 8.5\text{cm} \times 7.5\text{cm}$ (width, length, height) and these can move within a rectangular {\em work area} of $3\text{m} \times 2\text{m}$.  The platform is equipped with cameras with a top view to track the robots' positions.  Further, the Robotarium allows users to consider for control design a single integrator dynamics rather than the unicycle dynamics (the platform provides support functions to map single integrator dynamics to the unicycle dynamics). Hence,  the dynamics of the robot we want to control is $\bv{x}_{k} = \bv{x}_{k-1}+\bv{u}_{k}dt$,  where $\bv{x}_{k}=[p_{x,k},p_{y,k}]^{T}$ is the position of the robot at time-step $k$,  $\bv{u}_k=[v_{x,k},v_{y,k}]^{T}$ is the input vector of velocities (in the platform,  each  of the components are constrained to have modulus less than $0.5$m/s)  and $dt=0.033s$ is the Robotarium time-step.  

\begin{figure*}[thbp]
    \centering
    \includegraphics[width=0.62\columnwidth]{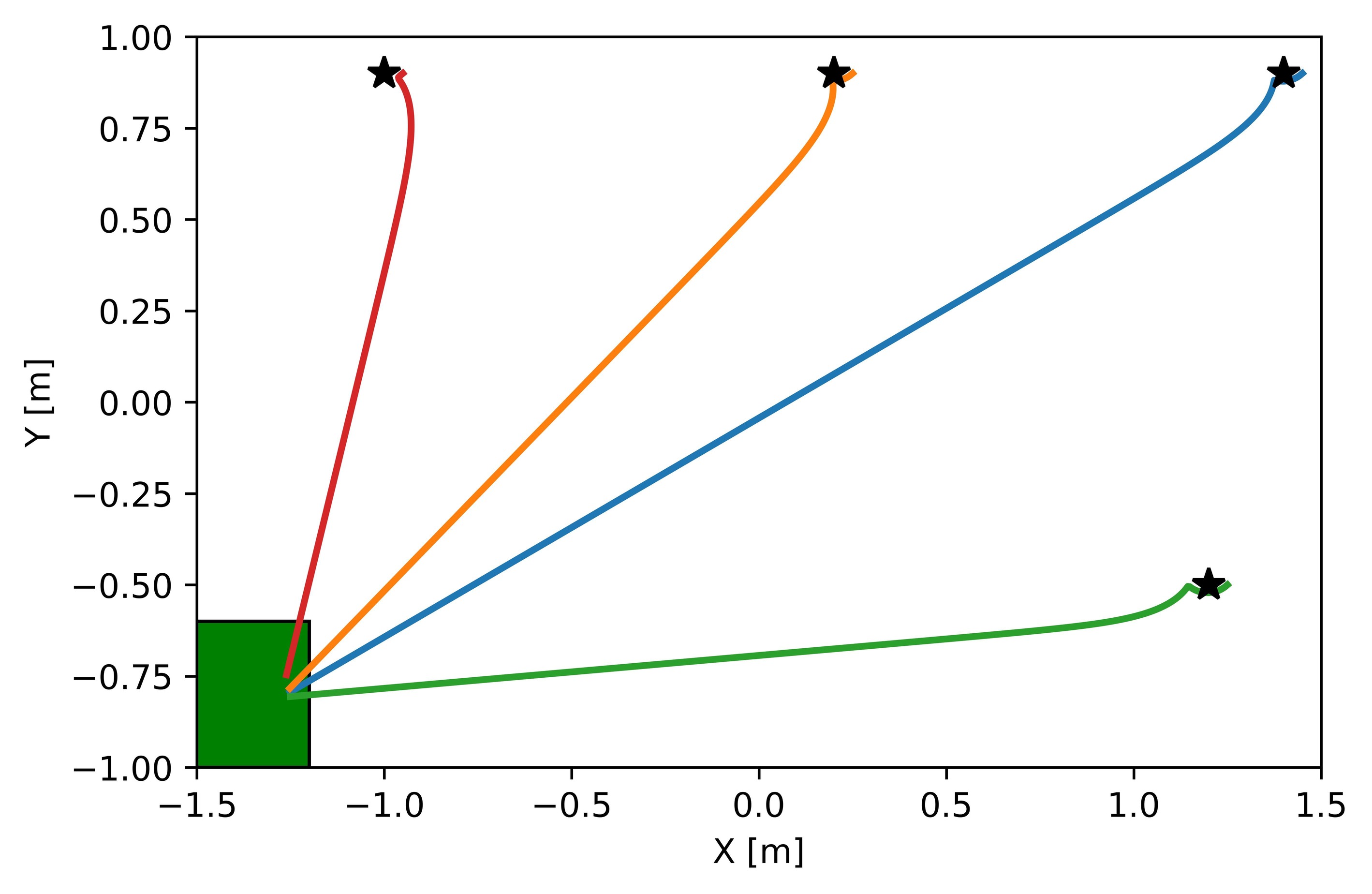} 
        \includegraphics[width=0.6\columnwidth]{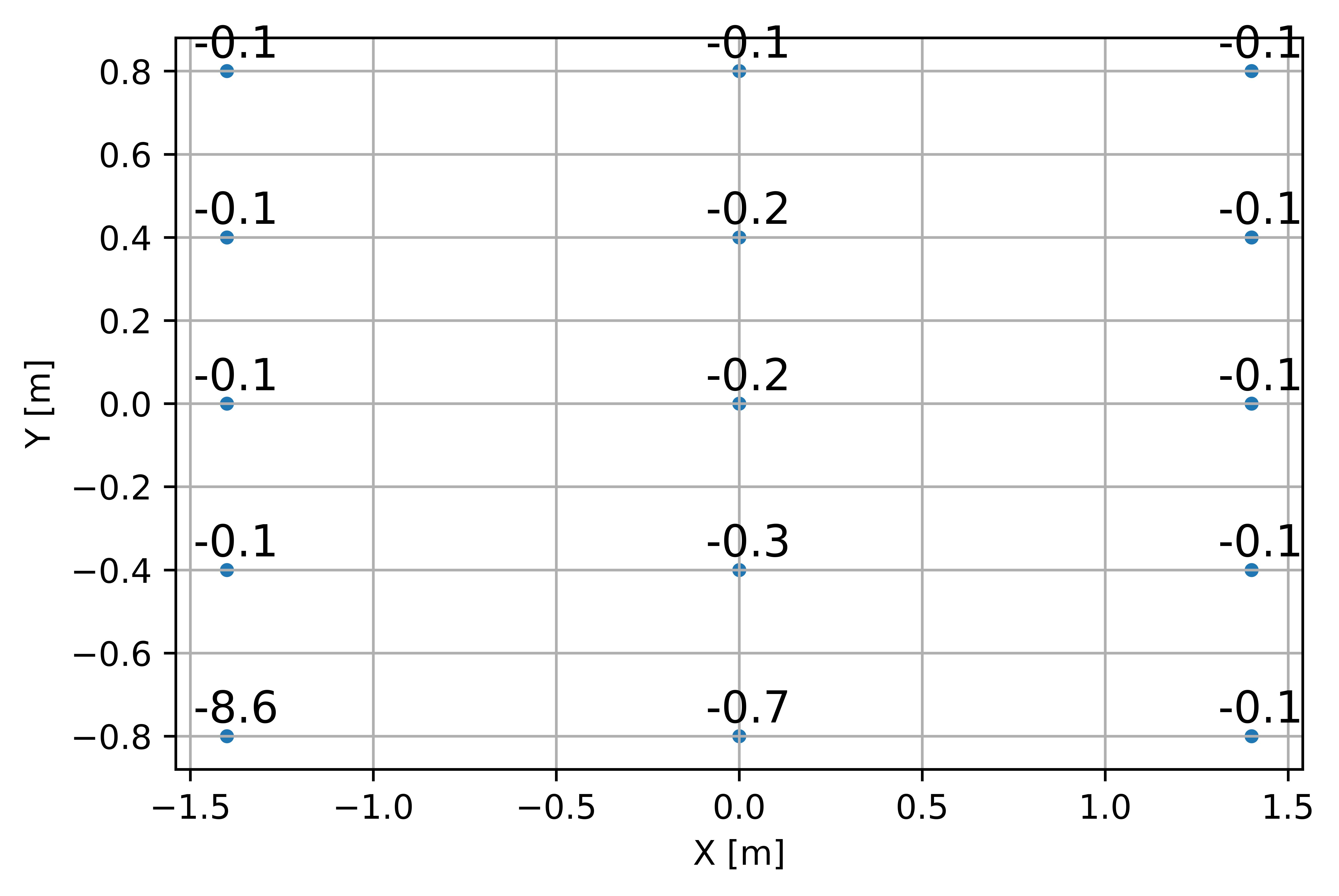} \\
         \includegraphics[width=0.67\columnwidth]{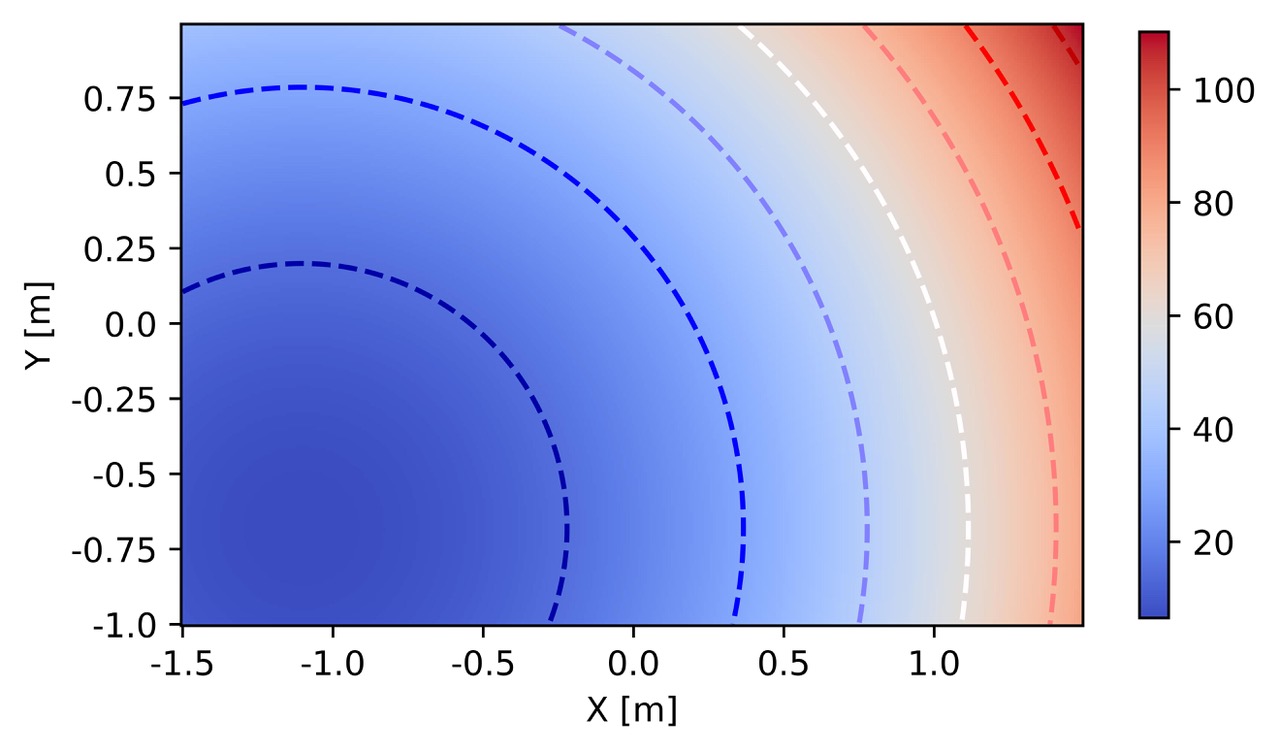} 
        \includegraphics[width=0.6\columnwidth]{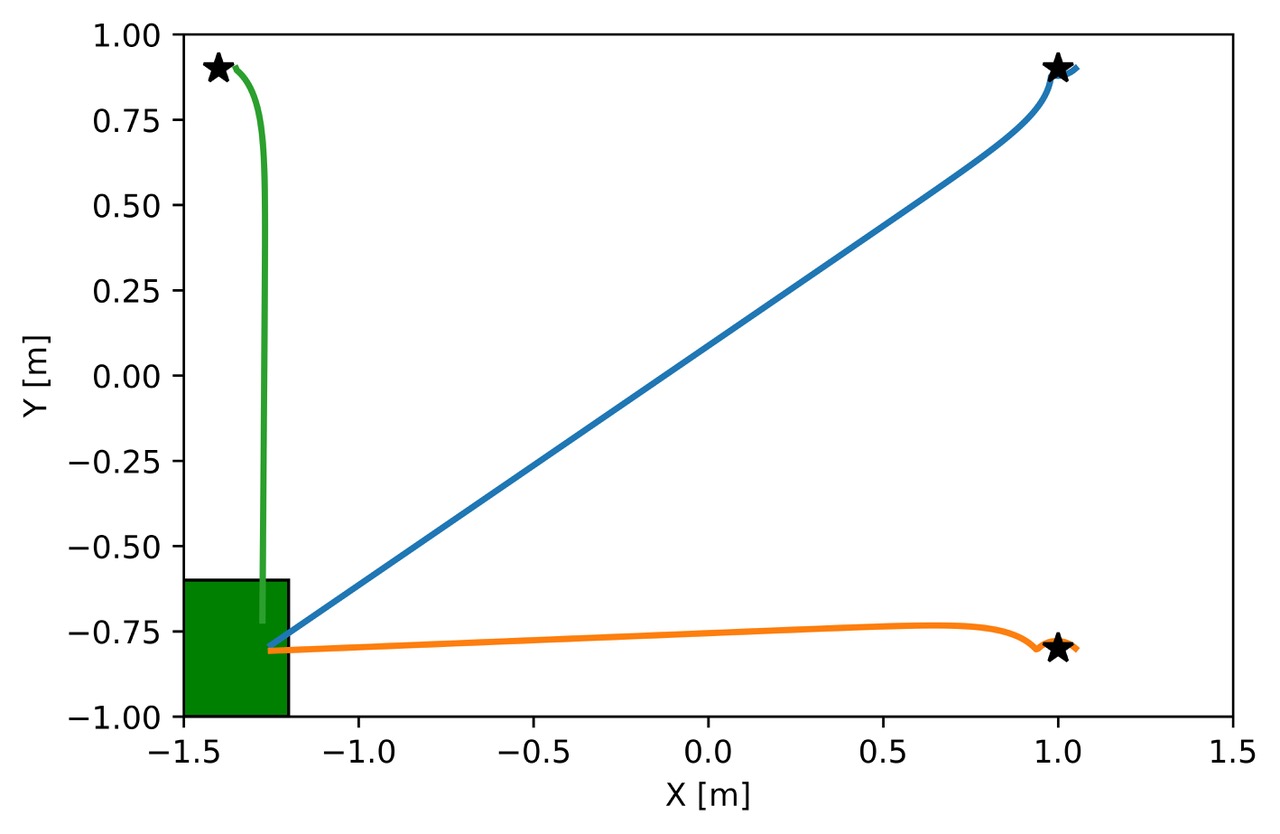} 
        \caption{{Top-left: robot trajectories starting from different initial positions ($\star$) when the policy in \eqref{eqn:gaussian_policy} - \eqref{eqn:gaussian_policy_recursion} is used (with $N=1$).  Top-right: the $\bv{o}_i$'s  together with the weights obtained via Algorithm \ref{alg:estimator}. Bottom: reconstructed cost (left) and robot trajectories when Algorithm \ref{alg:main} is used with this cost. Robot starts from initial positions that are different from those in the top panel.}}
    \label{fig:robot_Trajectories_Cont}
\end{figure*}

\begin{figure*}[thbp]
    \centering
    \includegraphics[width=0.65\columnwidth]{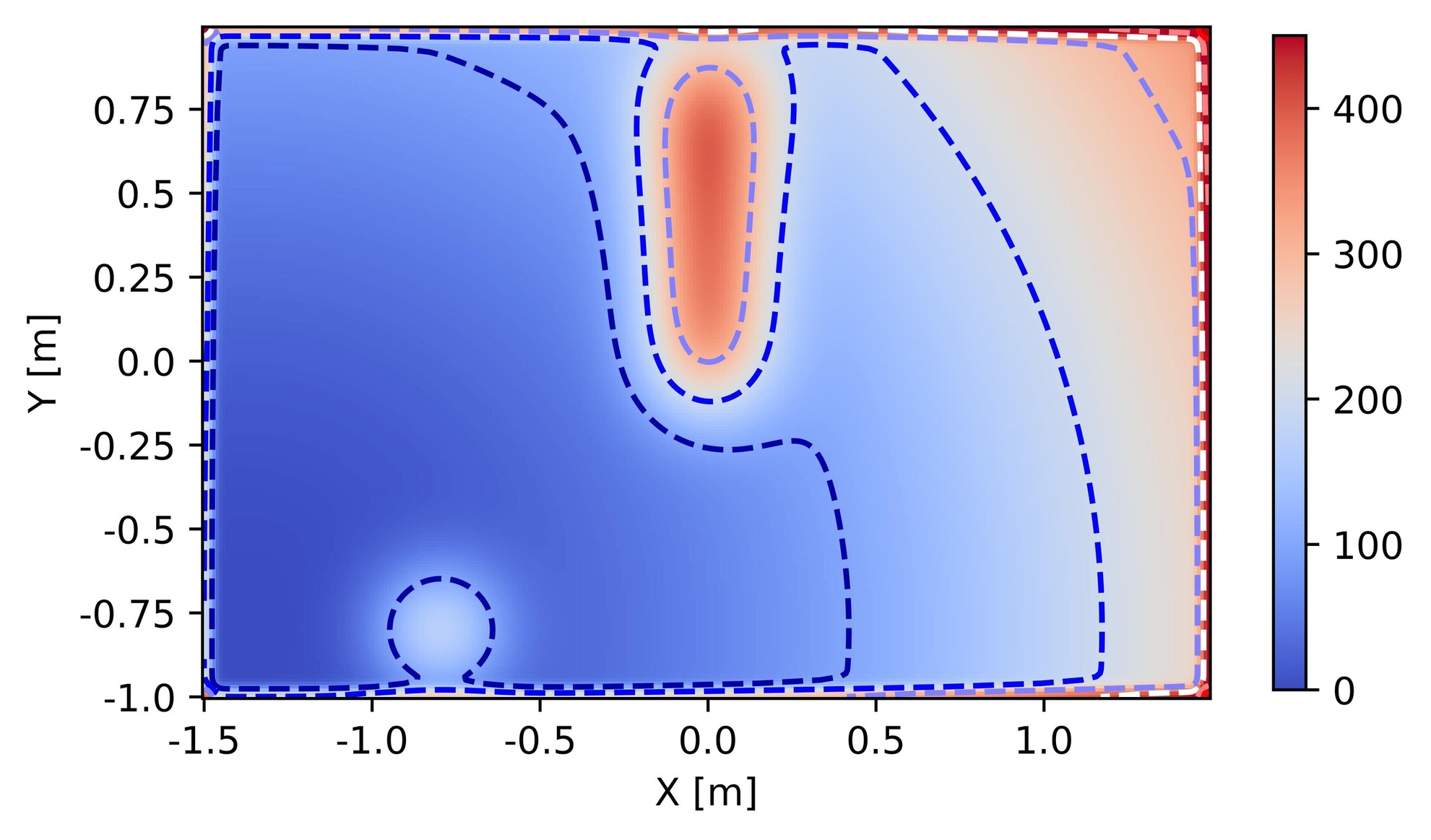} 
        \includegraphics[width=0.6\columnwidth]{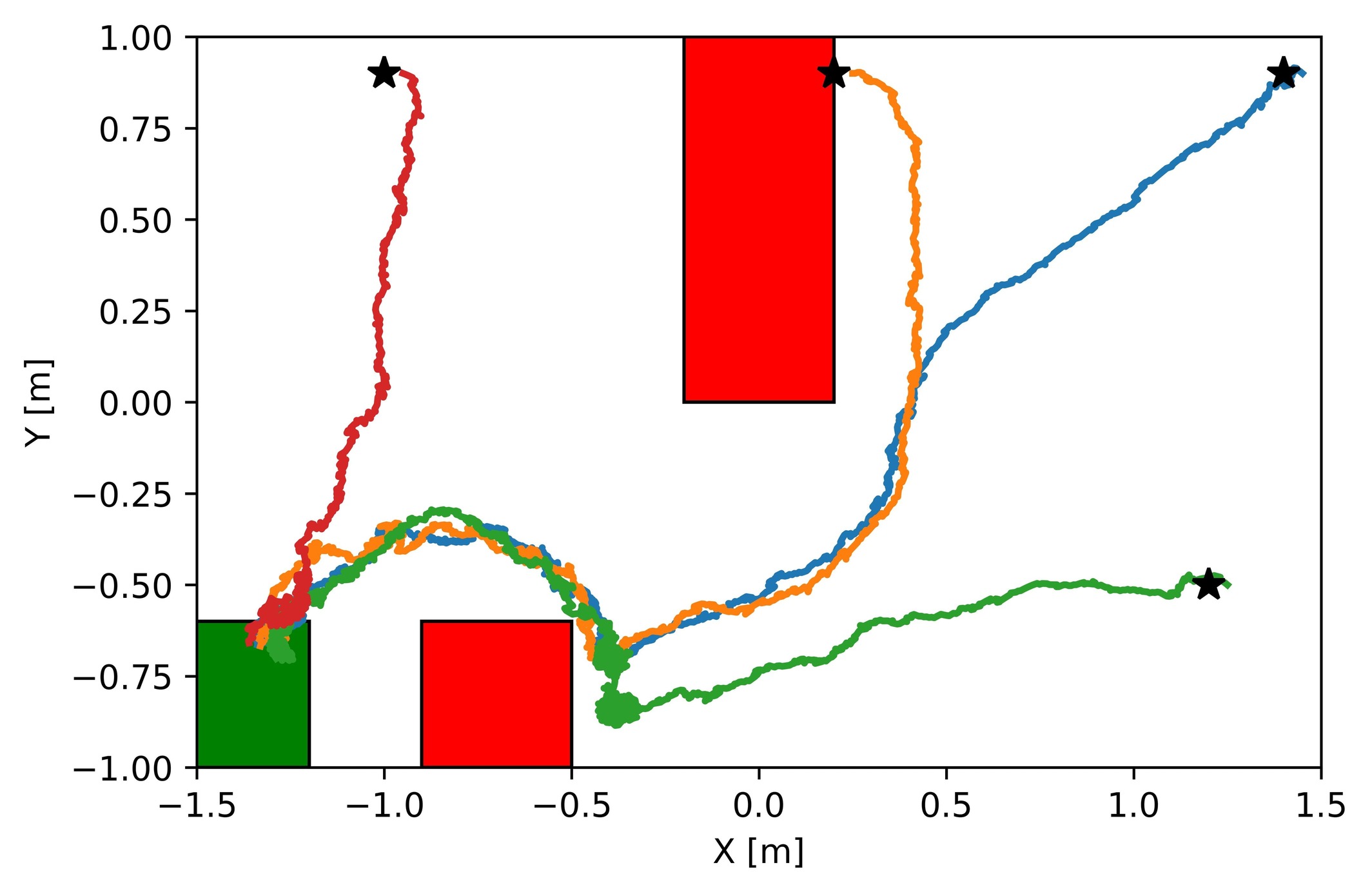} \\
        \includegraphics[width=0.65\columnwidth]{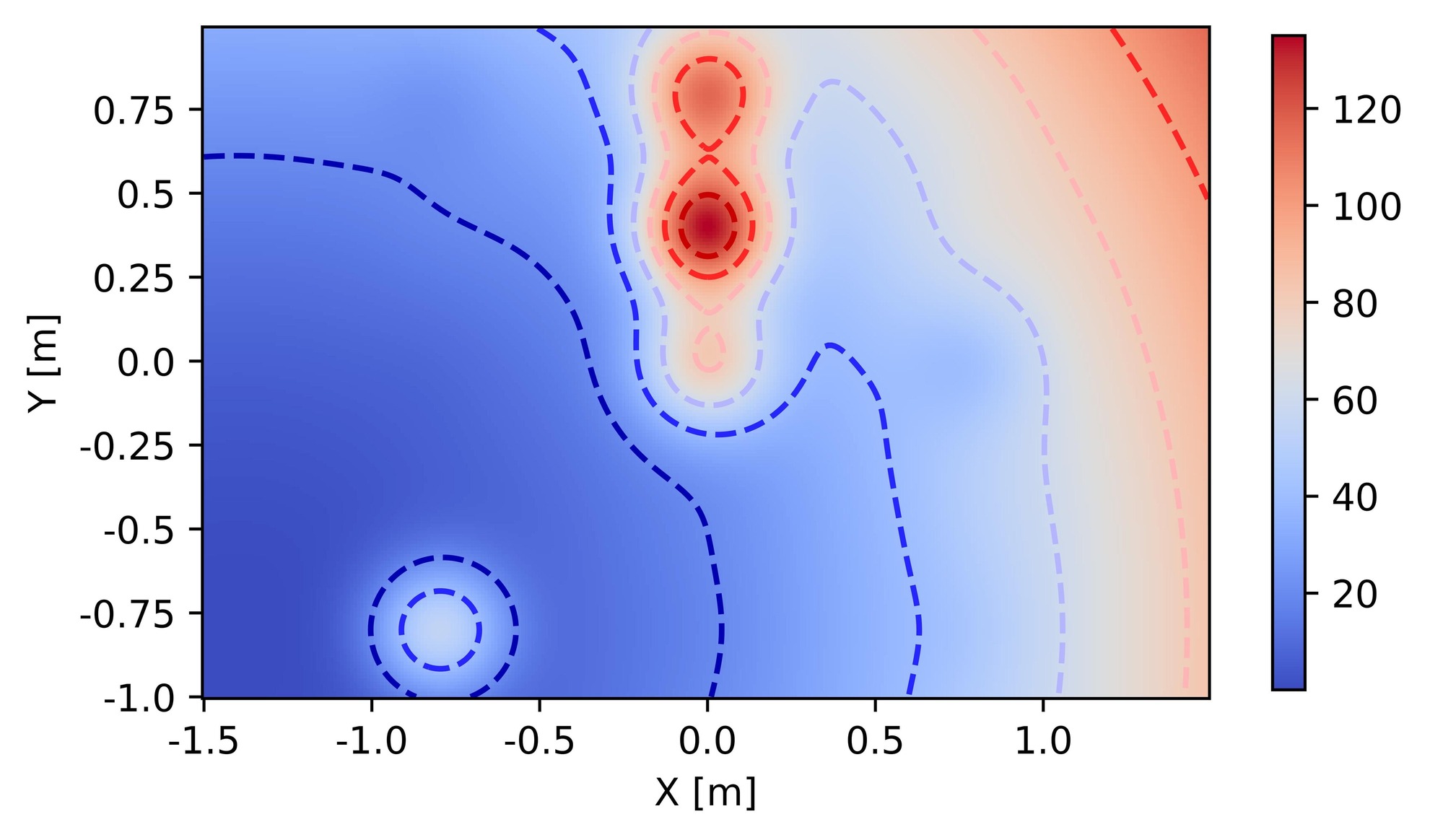} 
        \includegraphics[width=0.6\columnwidth]{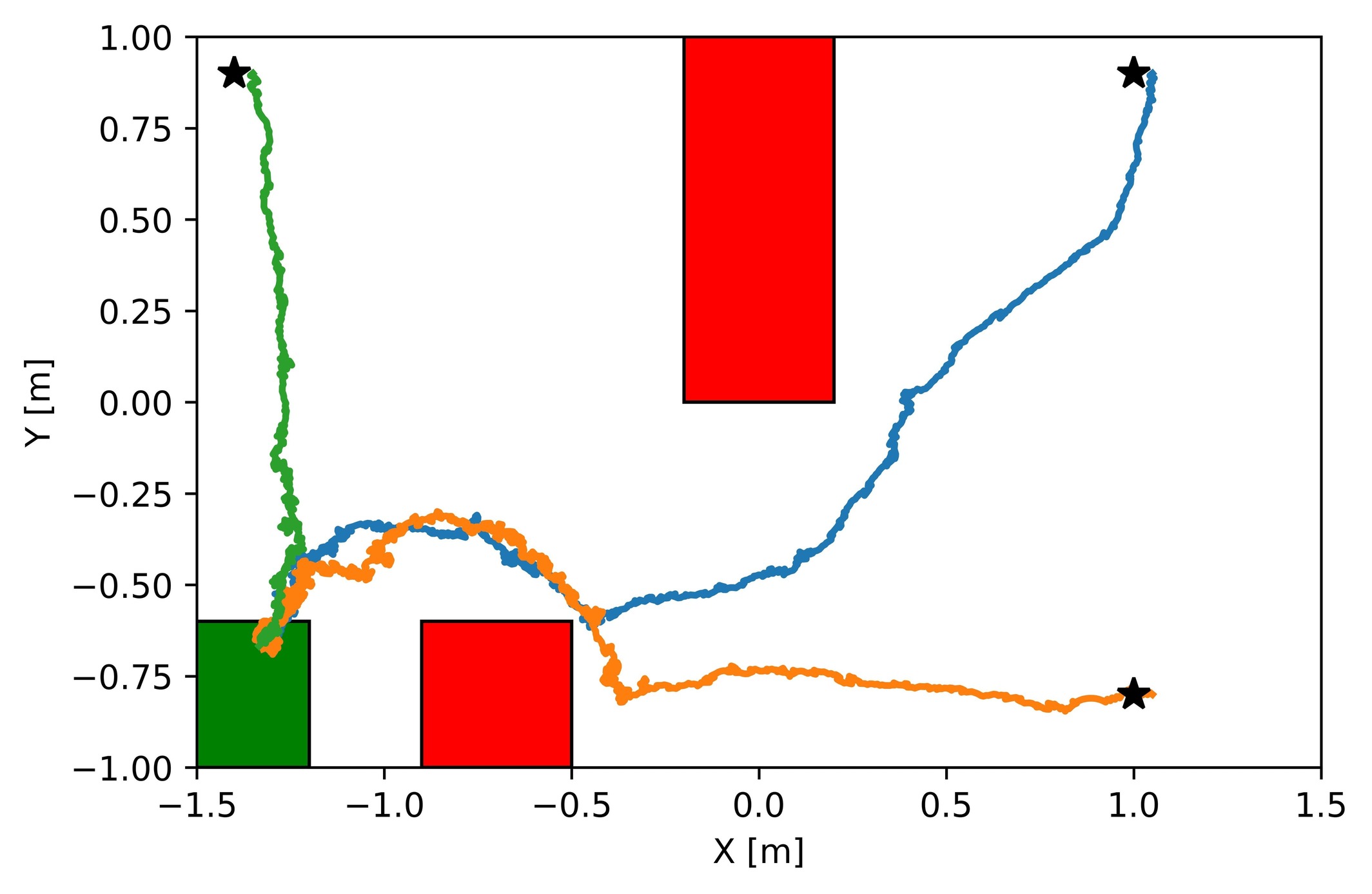} 
    \caption{Top-left: cost  for the {FOC} problem. Top-right: robot trajectories when the policy from Algorithm \ref{alg:main} is used  {(same initial positions and destination of Scenario $1$)}. Bottom panels: cost reconstructed via Algorithm \ref{alg:estimator} (left) and robot trajectories when Algorithm \ref{alg:main} is used with the estimated cost. Robots start from initial positions that are different from these in the top panel.}
    \label{fig:Robotarium_Results}
\end{figure*}

For the application of our results we set $\bv{x}_{k}\in\mathcal{X}:=[-1.5,1.5]\times[-1,1]$ and $\bv{u}_{k}\in\mathcal{U}:=[-0.5,0.5]\times[-0.5,0.5]$. Moreover, in our experiments we emulated measurement noise for the robot position that we assumed to be Gaussian as in e.g.,  \cite{7743685}. Thus,  $\plant{k|k-1}$ is given by $\mathcal{N}\left(\bv{x}_{k-1}+\bv{u}_{k}dt,\Sigma\right)$,  where we set 
$\bv{\Sigma} = \begin{bmatrix}
 0.001 & 0.0002  \\
 0.0002 & 0.001 \\
\end{bmatrix}$.\\
\\
{Given this set-up, we discuss separately the scenarios.\\
\\
\noindent {\bf Scenario $1$.} The desired destination for the robot, where it needs to stop, is $\bv{x}_d = [-1.4,-0.9]^T$, in the bottom-left corner of the Robotarium work area. Consequently,  for the FOC we set $c(\bv{x}_{k}) = 0.5(\bv{x}_k-\bv{x}_d)^T(\bv{x}_k-\bv{x}_d)$, $q^{(x)}_{k\mid k-1} = \mathcal{N}(\bv{x}_d,\bv{\Sigma}_x)$ and $q^{(u)}_{k\mid k-1} = \mathcal{N}(\bv{u}_d,\bv{\Sigma}_u)$,  where $\bv{u}_d = [0,0]^T$,  $\bv{\Sigma}_x = 0.003\cdot \bv{I}_{2}$ and $\bv{\Sigma}_u = 0.007\cdot\bv{I}_{2}$ ($\bv{I}_2$ is the $2\times 2$ identity matrix).  We then used the policy in \eqref{eqn:gaussian_policy} - \eqref{eqn:gaussian_policy_recursion} to control the robot and the results are in Figure \ref{fig:robot_Trajectories_Cont} (top-left panel).  In the figure, the behavior of the robot controlled with our policy is shown  when this starts from $4$ different positions. Next, we used the data pairs collected from the $4$ experiments to reconstruct the cost.  We used Algorithm \ref{alg:estimator} with a $15$-dimensional feature vector  ($F = 15$) and the $h_i$'s given by $h_i(\bv{x}_k) = (\bv{x}_k-\bv{o}_i)^T(\bv{x}_k-\bv{o}_i)$, where the $\bv{o}_i$'s were the coordinates of the points in the  work area shown in Figure \ref{fig:robot_Trajectories_Cont} (top-right).  The figure also shows the weights  returned by Algorithm \ref{alg:estimator} for each of the features.  The weights were obtained by adding a non-positivity constraint on $\bv{w}_{s}$ in the formulation of  \eqref{eqn:estimation_stationary}, leveraging the fact that the optimization problem (which we solved via CVX) would still be convex (see Remark \ref{cons_wi}).  The reconstructed cost is also shown in the bottom-left panel of Figure \ref{fig:robot_Trajectories_Cont}. We then gave back this cost to Algorithm \ref{alg:main} to verify if the robot would still be able to achieve the destination.  We performed these experiments by letting the robot start in positions that were different from the ones contained in the dataset used to estimate the cost. The results in Figure \ref{fig:robot_Trajectories_Cont} (bottom-right) confirm successful routing of the robot.  See our github for the documented code associated to the experiments described above. }\\
~
\\
{\noindent {\bf Scenario $2$}}. We decided not to build a target pf and  set $\idealplant{k|k-1}$ and $\idealcontrol{k|k-1}$ to be  uniform. Further, we discretized $\mathcal{U}$ in a $5\times5$ grid.  For the forward control problem, we used $c(\bv{x}_{k}) = 30(\bv{x}_{k} - \bv{x}_{d})^{2}+20\sum_{i=1}^{n}g_{i}(\bv{x}_{k}){+10b(\bv{x}_{k})}$,
where: (i) $\bv{x}_{d}$ is the desired goal/destination for the robot so that the first term in the cost promotes reaching the  goal; (ii) $n$ is the number of obstacles; (iii) the terms in the sum promote obstacle avoidance; {(iv) the last term is a cost associated to the boundaries of the Robotarium work area (the specific expression of this component of the cost is given on our github).}  We picked the $g_i$'s as
\begin{equation}\label{eqn:Gaussians_features}
 \begin{split}
&g_{i}(\bv{x}_{k}) := 
\frac{\exp\left(-0.5(\bv{x}_{k}-\bv{o}_{i})^T \bv{\Sigma}_o^{-1} (\bv{x}_{k}-\bv{o}_{i})\right)}{\sqrt{{(2\pi)^{2} \det(\bv{\Sigma}_o)}}} ,
 \end{split}
\end{equation}
with $\mathbf{o}_{i}$ being the coordinates of the barycenter of the $i$-{th} obstacle and $\bv{\Sigma}_o$ being the co-variance matrix.  See our github for the specific values of these parameters.  A plot of the cost function is given in Figure \ref{fig:Robotarium_Results} (top-left panel). We then used Algorithm \ref{alg:main} with $N=1$ to control the robot.  In our implementation of the algorithm, available at our github, the integrals to compute the expectations in \eqref{eqn:optimal_solution_statement_uniform} were estimated via Monte Carlo sampling and, since $\plant{k|k-1}$ is Gaussian, we used the analytic expression for the entropy.  We ran $4$ experiments from different initial robot positions {(the same used in Scenario $1$)} and,  for each of the experiments we recorded the robot trajectories.  In all the experiments,  as shown in  the top-right panel of Figure \ref{fig:Robotarium_Results},  the algorithm properly routed the robot to the destination while avoiding the obstacles. A video recording from an experiment is also given on our github. 
~
Next, we leveraged Algorithm \ref{alg:estimator} to reconstruct the cost using the data {pairs} collected from {the $4$ experiments above}.  We used a $16$-dimensional features vector (i.e., $F=16$), with the first $h_i$  equal to $(\bv{x}_{k}-\bv{x}_{d})^{2}$ and with the other $h_i$'s being Gaussians of the form \eqref{eqn:Gaussians_features} but centered around  the $15$ uniformly distributed points in the work area {of Figure \ref{fig:robot_Trajectories_Cont} (top-right).}  Again, we used CVX to solve the underlying optimization problem and, conveniently, this returned a vector of weights that would assign high values in magnitude to weights that corresponded to the obstacles and to the final destination. See the github for the values.  Finally,  we {again} gave back the reconstructed cost as an input to Algorithm \ref{alg:main}  with the aim of verifying if the robot would still be able to achieve the destination{, this time also} avoiding obstacles. {As for Scenario $1$}, we performed these experiments by letting the robots start {from new initial} positions  {and t}he results in Figure \ref{fig:Robotarium_Results} confirm successful routing and obstacle avoidance.  See \url{https://tinyurl.com/46uccxsf} for a recording from one  experiment.

\section{Conclusions}\label{sec:conclusions}
We considered the problem of reconstructing the possibly non-convex and non-stationary cost driving the actions of an agent from observations of its interactions with a nonlinear, non-stationary, and stochastic environment.  By leveraging  probabilistic descriptions that can be derived both directly from the data and from first-principles, we presented a result to tackle the inverse problem by solving a convex optimization problem.  To obtain this result we also {studied} a forward control problem with randomized policies as decision variables.  The results were also turned into algorithmic procedures with the code for our experiments made openly available. The effectiveness of our algorithms was experimentally evaluated via in-silico and hardware validations.  As part of our future work,  
{we aim to: (i) combine the  reproducing kernel Hilbert space and multiple kernel learning frameworks with our approach; (ii) use alternative {cost} parametrizations  that can be tackled via e.g.,  distributed solution methods \cite{9139372}{; (iii) study constrained versions of  Problem \ref{prob:main}}.  We also seek to: (i) engineer a principled feedback loop where the outcome of the inverse problem promotes exploration in the forward problem \cite{NIPS2012_2a50e9c2}; (ii) {following our discussion in Section \ref{sec:main_problem},} formalize, and exploit, links between our problem formulations, {fully probabilistic design} and optimal transport \cite{terpin2023dynamic}; (iii) extend our results to multi-agent systems performing repetitive tasks \cite{9993319,9317713} and to online settings \cite{6716965}.} From the applications viewpoint, we {aim to} design incentive schemes in sharing economy settings \cite{Sharing_book}.  

\section*{Acknowledgments}
EG and GR would like to thank Dr. Guy and Prof.  K\'arn\'y (Institute of Information Theory and Automation at the Czech Academy of Sciences) for the insightful discussions on a preliminary version of the results presented here. GR also wishes to thank Prof. Chiuso (Univ. of Padova) for the interesting insights on guided exploration. {The authors also wish to thank the AE and two anonymous reviewers who made several helpful comments and suggestions. These inputs led to improvements over the originally submitted manuscript.}

\appendix
\section{Appendix}
\subsection{Proof of Theorem \ref{thm:prob_main}}
\noindent{\bf Step $1$.} By means of Lemma \ref{lem:splitting_property}, Problem \ref{prob:main} can be recast as the sum of the following two sub-problems:
\begin{subequations}
\begin{equation}\label{eqn:split_1_problem_1}
    \begin{aligned}
    \underset{\{\control{k\mid k-1}\}_{1:N-1}}{\text{min}}
    &\bigg\{\DKL\left(p_{0:N-1}\mid\mid q_{0:N-1}\right) \\ 
    &  +\sum_{k=1}^{N-1} \E_{\bar{p}_{k-1:k}}\left[\E_{p^{(x)}_{k\mid k-1}}\left[c_k(\bv{X}_k)\right]\right]\bigg\} \\
   s.t. & \ \control{k\mid k-1}\in\sD \ \ \forall k\in 1:N-1,
    \end{aligned} 
\end{equation}
\text{and}
\begin{equation}\label{eqn:split_1_problem_2}
    \begin{aligned}
    \underset{\control{N\mid N-1}}{\text{min}}
    &\Big\{\E_{\bar{p}_{N-1}}\left[\DKL\left(p_{N\mid N-1}\mid\mid q_{N\mid N-1}\right)\right.\\ 
    &\left.+ \E_{p_{N\mid N-1}}\left[c_N(\bv{X}_N)\right]\right] \Big\}\\
   s.t. & \ \control{N\mid N-1}\in\sD.
    \end{aligned} 
\end{equation}
\end{subequations}  
Hence,  to solve Problem \ref{prob:main} we  solve (\ref{eqn:split_1_problem_2}) and take its solution into account in (\ref{eqn:split_1_problem_1}).  We let $\sC_N\left\{\control{N\mid N-1}\right\} :=\DKL\left(p_{N\mid N-1}\mid\mid q_{N\mid N-1}\right) +  \E_{p_{N\mid N-1}}\left[c_N(\bv{X}_N)\right]$.
Then,  the minimum of (\ref{eqn:split_1_problem_2}) is $\E_{\bar{p}_{N-1}}\left[\sC_N\left\{\optimalcontrol{N\mid N-1}\right\}\right],$ with 
\begin{equation}\label{eqn:split_1_problem_2_revised}
 \begin{aligned}
      \optimalcontrol{N\mid N-1} \in \underset{\control{N\mid N-1}}{\text{arg min}}
    &\Big\{\DKL\left(p_{N\mid N-1}\mid\mid q_{N\mid N-1}\right)\\ 
    & +  \E_{p_{N\mid N-1}}\left[\bar{c}_N(\bv{X}_N)\right]\Big\} \\
   s.t. & \ \control{N\mid N-1}\in\sD,
    \end{aligned} 
\end{equation}
 where we set $\bar{c}_N(\bv{x}_N):={c}_N(\bv{x}_N) + \hat{c}_N(\bv{x}_N)$, $\hat{c}_N(\bv{x}_N)=0$. This corresponds to the recursion in (\ref{eqn:backward_recursion}) at $k=N$. 

\noindent{\bf Step $2$.}The constraint in (\ref{eqn:split_1_problem_2_revised}) is linear in the decision variable and we show convexity of the problem by showing that its cost functional is convex.  {N}ote that:
\begin{align*}
&\sC_N\{\control{N\mid N-1}\} {=}\\
&\DKL(p_{N\mid N-1}\|q_{N\mid N-1}) + \E_{p_{N\mid N-1}}[\bar{c}_N(\bv{X}_N)]\\
& =  \E_{\control{N\mid N-1}}\Big[\DKL\left(\plant{N\mid N-1}\|\idealplant{N\mid N-1}\right)\\
& +\E_{\plant{N\mid N-1}}[\bar{c}_N(\bv{X}_N)]\Big] + \DKL\left(\control{N\mid N-1}\|\idealcontrol{N\mid N-1}\right).
\end{align*}
Hence, the problem in (\ref{eqn:split_1_problem_2_revised}) becomes
\begin{equation}\label{eqn:split_1_problem_2_revised_for_convex}
\begin{aligned}
\underset{\control{N\mid N-1}}{\text{min}}
& \bigg\{\E_{\control{N\mid N-1}}\Big[\DKL\left(\plant{N\mid N-1} \mid \mid \idealplant{N\mid N-1}\right) \\
&+ \E_{\plant{N\mid N-1}}\left[\bar{c}_N(\bv{X}_N)\right]\Big] \\
& + \DKL\left(\control{N\mid N-1} \mid \mid \idealcontrol{N\mid N-1}\right) \bigg\}\\
s.t. & \ \control{N\mid N-1} \in \sD,
\end{aligned}
\end{equation}
which has a strictly convex cost.

\noindent{\bf Step $3$.} 
The Lagrangian of the problem in (\ref{eqn:split_1_problem_2_revised_for_convex}) is $\sL(\control{N\mid  N-1},\lambda_N) = \E_{\control{N\mid  N-1}}\Big[\DKL\left(\plant{N\mid  N-1}\mid  \mid  \idealplant{N\mid  N-1}\right) +\E_{\plant{N\mid  N-1}}\left[\bar{c}_N(\bv{X}_N)\right]\Big] + \DKL\left(\control{N\mid  N-1}\mid  \mid  \idealcontrol{N\mid  N-1}\right)  +\lambda_N\left(\sum_{\bv{u}_k}\control{N\mid  N-1}-1\right)$, where $\lambda_N$ is the Lagrange multiplier corresponding to the constraint $\control{N\mid  N-1}\in\sD$. We find the optimal solution by imposing the first order stationarity conditions on $ \sL(\control{N\mid  N-1},\lambda_N)$. We start with imposing the stationarity condition with respect to $\control{N\mid  N-1}$. This yields $\DKL\left(\plant{N\mid  N-1}\mid  \mid  \idealplant{N\mid  N-1}\right)+\E_{\plant{N\mid  N-1}}\left[\bar{c}_N(\bv{X}_N)\right] +\ln\frac{\control{N\mid  N-1}}{\idealcontrol{N\mid  N-1}}+1+\lambda_N = 0$.
Hence, any candidate solution is of the form:
\begin{equation}\label{eqn:candidate_solution}
\candidatecontrol{N\mid  N-1} = \frac{\bar{p}_{N\mid  N-1}^{(u)}\exp\left(-\E_{\plant{N\mid  N-1}}\left[\bar{c}_N(\bv{X}_N)\right]\right)}{\exp\left(1+\lambda_N\right)},
\end{equation}
where $\bar{p}_{N\mid  N-1}^{(u)}$ is defined as in \eqref{eqn:p_bar}.
Now, by imposing the stationarity condition for $ \sL(\control{N\mid  N-1},\lambda_N)$ with respect to $\lambda_N$, we get $\sum_{\bv{u}_N}\control{N\mid  N-1}-1 = 0$.  Hence, we get:
\begin{equation}\label{eqn:stationarity_cond_Lagrange_multiplier}
\begin{split}
&\sum_{\bv{u}_N}\bar{p}_{N\mid  N-1}^{(u)}\exp\Big(-\E_{\plant{N\mid  N-1}}\left[\bar{c}_N(\bv{X}_N)\right]\Big)= \exp\left(1+\lambda_N\right).
\end{split}
\end{equation}
Since the problem in (\ref{eqn:split_1_problem_2_revised_for_convex}) is convex with a strictly convex cost functional, (\ref{eqn:candidate_solution}) and (\ref{eqn:stationarity_cond_Lagrange_multiplier}) imply \cite{Roc_88} that the unique optimal solution  is
\begin{align}\label{eqn:optimal_solution_N}
&\optimalcontrol{N\mid  N-1}  = \frac{\bar{p}_{N\mid  N-1}^{(u)}\exp\left(-\E_{\plant{N\mid  N-1}}\left[\bar{c}_N(\bv{X}_N)\right]\right)}{\sum_{\bv{u}_N}\bar{p}_{N\mid  N-1}^{(u)}\exp\left(-\E_{\plant{N\mid  N-1}}\left[\bar{c}_N(\bv{X}_N)\right]\right)}.
\end{align}
This is  the optimal solution given in (\ref{eqn:optimal_solution_statement}) for $k=N$, with $\bar{c}_N(\bv{x}_N)$ generated via the backward recursion in (\ref{eqn:backward_recursion}). Moreover,  the minimum of the problem in (\ref{eqn:split_1_problem_2_revised_for_convex}) is:
\begin{equation*}
\begin{split}
& \sC_N\left\{\optimalcontrol{N\mid  N-1}\right\}  =\\
&-\ln\Bigg(\sum_{\bv{u}_N}\idealcontrol{N\mid  N-1}\exp\Big(-\DKL\left(\plant{N\mid  N-1}\mid  \mid  \idealplant{N\mid  N-1}\right)\\
&-\E_{\plant{N\mid  N-1}}\left[\bar{c}_N(\bv{X}_N)\right]\Big)\Bigg).
\end{split}
\end{equation*}
Hence, the minimum for the sub-problem in (\ref{eqn:split_1_problem_2}) is
\begin{equation}\label{eqn:cost_N}
\begin{split}
& \E_{\bar{p}_{N-1}}\left[\sC_N\left\{\optimalcontrol{N\mid  N-1}\right\}\right]  = -\E_{\bar{p}_{N-1}}\left[\hat{c}_{N-1}(\bv{X}_{N-1})\right],
\end{split}
\end{equation}
where 
\begin{equation*}
\begin{split}
    &\hat{c}_{N-1}(\bv{x}_{N-1}):=\\ 
    &\ln\Big(\E_{\idealcontrol{N\mid  N-1}}\Big[\exp\Big(-\DKL\left(\plant{N\mid  N-1}\mid  \mid  \idealplant{N\mid  N-1}\right)\\
    &-\E_{\plant{N\mid  N-1}}\left[\bar{c}_N(\bv{X}_N)\right]\Big)\Big]\Big).
\end{split}
\end{equation*}
This is the optimal cost for $k=N$ given in part (ii) of the statement. Next, we make use of the minimum found for the sub-problem (\ref{eqn:split_1_problem_2}) to solve the sub-problem corresponding to $k\in 1:N-1$.

\noindent{\bf Step $4$.} Since the problem in (\ref{eqn:main_problem}) has been split as the sum of the sub-problems in (\ref{eqn:split_1_problem_1}) - (\ref{eqn:split_1_problem_2}) and since the solution of (\ref{eqn:split_1_problem_2}) gives the minimum (\ref{eqn:cost_N}), we have that Problem \ref{prob:main} is equivalent to
\begin{equation}\label{eqn:problem_N-1}
    \begin{aligned}
    \underset{\{\control{k\mid  k-1}\}_{1:N-1}}{\text{min}}
    & \Bigg\{\DKL\left(p_{0:N-1}\mid  \mid  q_{0:N-1}\right)\\  
    &+\sum_{k=1}^{N-1} \E_{\bar{p}_{k-1:k}}\left[\E_{p^{(x)}_{k\mid  k-1}}\left[c_k(\bv{X}_k)\right]\right]\\  
    &-\E_{\bar{p}_{N-1}}\left[\hat{c}_{N-1}(\bv{X}_{N-1})\right]\Bigg\}  \\
   s.t. & \ \control{k\mid  k-1}\in\sD \ \ \forall k\in 1:N-1.
    \end{aligned} 
\end{equation}
In turn,  the cost of the above problem can be written as:
\begin{equation}\label{eqn:decomposition_N-1}
    \begin{split}
&\DKL\left(p_{0:N-2}\mid  \mid  q_{0:N-2}\right) +\sum_{k=1}^{N-2} \E_{\bar{p}_{k-1:k}}\left[\E_{p^{(x)}_{k\mid  k-1}}\left[c_k(\bv{X}_k)\right]\right]\\ 
&+ \E_{p_{0:N-2}}\left[\DKL\left(p_{N-1\mid  N-2}\mid  \mid  q_{N-1\mid  N-2}\right) \right] \\
& +  \E_{\bar{p}_{N-2:N-1}}\left[\E_{p^{(x)}_{N-1\mid  N-2}}\left[c_{N-1}(\bv{X}_{N-1})\right]\right]\\
&-\E_{\bar{p}_{N-1}}\left[\hat{c}_{N-1}(\bv{X}_{N-1})\right],
    \end{split} 
\end{equation}
{where we recall from Section \ref{sec:set-up} that $p^{(x)}_{N-1\mid  N-2}$ is the shorthand notation for $p^{(x)}_{N-1}\left(\bv{x}_{N-1}\mid\bv{x}_{N-2},\bv{u}_{N-1}\right)$.  Also,  following e.g., \citep{gagliardi2020probabilistic}, for the last three terms in \eqref{eqn:decomposition_N-1} we note that: (i) $\DKL\left(p_{N-1\mid  N-2}\mid  \mid  q_{N-1\mid  N-2}\right)$ only depends on $\bv{X}_{N-2}$ and hence $
\E_{p_{0:N-2}}\left[\DKL\left(p_{N-1\mid  N-2}\mid  \mid  q_{N-1\mid  N-2}\right) \right]  = \E_{\bar{p}_{N-2}}\left[\DKL\left(p_{N-1\mid  N-2}\mid  \mid  q_{N-1\mid  N-2}\right) \right]$; (ii) the expectation on the third line of \eqref{eqn:decomposition_N-1} can be written as $\E_{\bar{p}_{N-2}}\left[\E_{p_{N-1\mid  N-2}}\left[c_{N-1}(\bv{X}_{N-1})\right]\right]$, where $p_{N-1\mid  N-2}$ is defined  (in Section \ref{sec:set-up}) as $p\left(\bv{x}_{N-1},\bv{u}_{N-1}\mid\bv{x}_{N-2}\right)$; (iii) the expectation in the last line of \eqref{eqn:decomposition_N-1}  can be recast as $\E_{\bar{p}_{N-2}}\left[\E_{p_{N-1\mid  N-2}}\left[\hat{c}_{N-1}(\bv{X}_{N-1})\right]\right]$. Therefore,}
the last three terms in (\ref{eqn:decomposition_N-1}) can be  written as 
\begin{equation*}
    \begin{aligned}
        &\E_{\bar{p}_{N-2}}\Big[\DKL\left(p_{N-1\mid  N-2}\mid  \mid  q_{N-1\mid  N-2}\right)\\
        & +  \E_{p_{N-1\mid  N-2}}\left[\bar{c}_{N-1}(\bv{X}_{N-1})\right]\Big],
    \end{aligned}
\end{equation*}
where $\bar{c}_{N-1}(\bv{x}_{N-1}):= c_{N-1}(\bv{x}_{N-1}) -\hat{c}_{N-1}(\bv{x}_{N-1})$. This corresponds to the recursion (\ref{eqn:backward_recursion}) at time $k=N-1$. Now, the problem in (\ref{eqn:problem_N-1}) can be again split, this time as the sum of the following two sub-problems:
\begin{subequations}
\begin{equation}\label{eqn:split_2_problem_1}
    \begin{aligned}
    \underset{\{\control{k\mid  k-1}\}_{1:N-2}}{\text{min}}
    &\Bigg\{\DKL\left(p_{0:N-2}\mid  \mid  q_{0:N-2}\right)\\ 
    & +\sum_{k=1}^{N-2} \E_{\bar{p}_{k-1:k}}\left[\E_{p^{(x)}_{k\mid  k-1}}\left[c_k(\bv{X}_k)\right]\right]\Bigg\}  \\
   s.t. & \ \control{k\mid  k-1}\in\sD \ \ \forall k\in 1:N-2,
    \end{aligned} 
\end{equation}
\text{and}
\begin{equation}\label{eqn:split_2_problem_2}
    \begin{aligned}
    \underset{\control{N-1\mid  N-2}}{\text{min}}
    &\Bigg\{\E_{\bar{p}_{N-2}}\big[\DKL\left(p_{N-1\mid  N-2}\mid  \mid  q_{N-1\mid  N-2}\right)\\
    &+  \E_{p_{N-1\mid  N-2}}\left[\bar{c}_{N-1}(\bv{X}_{N-1})\right]\big]\Bigg\} \\
   s.t. & \ \control{N-1\mid  N-2}\in\sD.
    \end{aligned} 
\end{equation}
\end{subequations} 
Again, the sub-problem in (\ref{eqn:split_2_problem_2}) is independent on the sub-problem in (\ref{eqn:split_2_problem_1}) and its decision variable, i.e. $\control{N-1\mid  N-2}$, is independent on $\bar{p}_{N-2}$. Hence the minimum is $\E_{\bar{p}_{N-2}}\left[\sC_{N-1}\left\{\optimalcontrol{N-1\mid  N-2}\right\}\right]$, with
\begin{equation*}
\begin{aligned}
    \sC_{N-1}\left\{\control{N-1\mid  N-2}\right\} & :=\DKL\left(p_{N-1\mid  N-2}|\mid  q_{N-1\mid  N-2}\right)\\
    &+  \E_{p_{N-1\mid  N-2}}\left[\bar{c}_{N-1}(\bv{X}_{N-1})\right],
\end{aligned}
\end{equation*}
and $\optimalcontrol{N-1\mid  N-2}$ being the solution of
\begin{equation}\label{eqn:split_2_problem_2_revised}
    \begin{aligned}
    \underset{\control{N-1\mid  N-2}}{\text{min}}
    &\Bigg\{\DKL\left(p_{N-1\mid  N-2}\mid  \mid  q_{N-1\mid  N-2}\right)\\
    & +  \E_{p_{N-1\mid  N-2}}\left[\bar{c}_{N-1}(\bv{X}_{N-1})\right]\Bigg\} \\
   s.t. & \ \control{N-1\mid  N-2}\in\sD.
    \end{aligned} 
\end{equation}
This has the same structure as  (\ref{eqn:split_1_problem_2_revised}) and hence the unique optimal solution for the sub-problem in (\ref{eqn:split_2_problem_2_revised}) is
\begin{multline}\label{eqn:optimal_solution_N-1}
\optimalcontrol{N-1\mid  N-2}=\\
\frac{\bar{p}_{N-1\mid  N-2}^{(u)}\exp\left(-\E_{\plant{N-1\mid  N-2}}\left[\bar{c}_{N-1}(\bv{X}_{N-1})\right]\right)}{\sum_{\bv{u}_{N-1}}\bar{p}_{N-1\mid  N-2}^{(u)}\exp\left(-\E_{\plant{N-1\mid  N-2}}\left[\bar{c}_{N-1}(\bv{X}_{N-1})\right]\right)}{.}
\end{multline}
That is, (\ref{eqn:optimal_solution_N-1}) is the optimal solution given in (\ref{eqn:optimal_solution_statement}) for $k=N-1$, with $\bar{c}_{N-1}(\bv{x}_{N-1})$ obtained via the backward recursion in (\ref{eqn:backward_recursion}). Moreover, the corresponding cost for the problem in (\ref{eqn:split_2_problem_2}) is $\E_{\bar{p}_{N-2}}\left[\sC_{N-1}\left\{\optimalcontrol{N-1\mid  N-2}\right\}\right]  = -\E_{\bar{p}_{N-2}}\left[\hat{c}_{N-2}(\bv{X}_{N-2})\right]$,
where
\begin{equation*}
\begin{split}
 &\hat{c}_{N-2}(\bv{x}_{n-2}) =\\ 
 & \ln\Big(\E_{\idealcontrol{N-1\mid  N-2}}\Big[\exp\Big(-\DKL\left(\plant{N-1\mid  N-2}\mid  \mid  \idealplant{N-1\mid  N-2}\right)\\
 &-\E_{\plant{N-1\mid  N-2}}\left[\bar{c}_{N-1}(\bv{X}_{N-1})\right]\Big)\Big]\Big).
 \end{split}
\end{equation*}
This is  the optimal cost for $k=N-1$ given in  part (ii) of the statement. We can now draw the desired conclusions.

\noindent{\bf Step $5$.} By iterating Step $4$, at each of the remaining time-steps in the window $1:N-2$, Problem \ref{prob:main} can always be split in sub-problems, where the sub-problem corresponding to the last time instant is given by:
 \begin{equation}\label{eqn:split_last_problem_general}
    \begin{aligned}
    \underset{\control{k\mid  k-1}}{\text{min}}&\E_{\bar{p}_{k-1}}\left[\DKL\left(p_{k\mid  k-1}\mid  \mid  q_{k\mid  k-1}\right) +  \E_{p_{k\mid  k-1}}\left[\bar{c}_{k}(\bv{X}_{k})\right]\right] \\
   s.t. & \ \control{k\mid  k-1}\in\sD,
    \end{aligned} 
\end{equation}
where $\bar{c}_{k}(\bv{x}_{k}):= c_{k}(\bv{x}_{k}) -\hat{c}_{k}(\bv{x}_{k})$ and,
\begin{equation*}
  \begin{aligned}
      \hat{c}_{k}(\bv{x}_{k})&=\ln\Bigg(\E_{\idealcontrol{k+1\mid  k}}\Big[\exp\Big(-\DKL\left(\plant{k+1\mid  k}\mid  \mid  \idealplant{k+1\mid  k}\right)\\
      & -\E_{\plant{k+1\mid  k}}\left[\bar{c}_{k+1}(\bv{X}_{k+1})\right]\Big)\Big]\Bigg).
  \end{aligned}  
\end{equation*}
This yields the recursion in (\ref{eqn:backward_recursion}) a time $k$. Hence, the optimal solution for the sub-problem is 
\begin{align*}
&\optimalcontrol{k\mid  k-1} = \frac{\bar{p}^{(u)}_{k\mid  k-1}\exp\left(-\E_{\plant{k\mid  k-1}}\left[\bar{c}_{k}(\bv{X}_{k})\right]\right)}{\sum_{\bv{u}_{k}}\bar{p}^{(u)}_{k\mid  k-1}\exp\left(-\E_{\plant{k\mid  k-1}}\left[\bar{c}_{k}(\bv{X}_{k})\right]\right)}.
\end{align*}
This is the optimal solution given in (\ref{eqn:optimal_solution_statement}) at time $k$, with $\bar{c}_{k}(\bv{x}_{k})$ obtained from the backward recursion in (\ref{eqn:backward_recursion}). Part (i) of the result is then proved.  The corresponding optimal cost at time $k$ is $-\E_{\bar{p}_{k-1}}\left[\hat{c}_{k-1}(\bv{X}_{k-1})\right]$. Thus, the optimal cost Problem \ref{prob:main} is $-\sum_{k=1}^N\E_{\bar{p}_{k-1}}\left[\hat{c}_{k-1}(\bv{X}_{k-1})\right]$ and this proves part (ii) of the result. \hspace*{\fill}\qed

{\subsection{Proofs of Instrumental Results}
\noindent{\bf Proposition \ref{prop:assumption}.} Item (ii) implies that the second part of the cost in \eqref{eqn:twisted_cost} is lower bounded by $-\bar{e}$ and upper bounded by $0$.   Hence, we need to show that (i) implies boundedness of the first component of \eqref{eqn:twisted_cost}.  By recursively applying Lemma \ref{lem:splitting_property}, we get that the first component of \eqref{eqn:twisted_cost} can be written as
$\DKL(p_0(\bv{x}_0)\mid\mid \tilde q_0(\bv{x}_0)) + \sum_{k=1}^N \E_{p_{k-1}(\bv{x}_{k-1})}\left[\DKL\left(p_{k\mid k-1}\mid\mid \tilde{q}_{k\mid k-1}\right)\right]$.  Now, (i) implies that $\DKL(p_0(\bv{x}_0)\mid\mid \tilde q_0(\bv{x}_0))\le H_0$. Additionally,  at each $k$, the term inside the expectation of the above expression can be written as $\E_{p^{(u)}_{k\mid k-1}}\left[\DKL\left(p^{(x)}_{k\mid k-1}\mid\mid \tilde{q}^{(x)}_{k\mid k-1}\right)\right] + \DKL\left(p^{(u)}_{k\mid k-1}\mid\mid \tilde{q}^{(u)}_{k\mid k-1}\right)$ and from  (i)  there exist some constant $H_k$ such that $\DKL\left(p^{(x)}_{k\mid k-1}\mid\mid \tilde{q}^{(x)}_{k\mid k-1}\right) \le H_k$. The desired conclusion then follows by noticing that, at each $k$,  the term $\DKL\left(p^{(u)}_{k\mid k-1}\mid\mid \tilde{q}^{(u)}_{k\mid k-1}\right) =0$ if $p^{(u)}_{k\mid k-1} = \tilde{q}^{(u)}_{k\mid k-1}$.
\\
\\
\noindent{\bf Derivations for the policy in \eqref{eqn:gaussian_policy} - \eqref{eqn:gaussian_policy_recursion}}.  We note that,  in the given setting,  $\DKL\left(\plant{k|k-1}||\idealplant{k|k-1}\right)= 0.5(\bv{A}\bv{x}_{k-1}+\bv{B}\bv{u}_{k}-\bv{x}_{d})^{T}\bv{R^{-1}}(\bv{A}\bv{x}_{k-1}+\bv{B}\bv{u}_{k}-\bv{x}_{d}) -0.5\Tr(\bv{\Sigma}(\bv{\Sigma}^{-1}-\bv{R}^{-1}))$ and $ \E_{\plant{k|k-1}}\left[c_{k}(\bv{x}_{k})\right] = 0.5(\bv{A}\bv{x}_{k-1}+\bv{B}\bv{u}_{k}-\bv{x}_{d})^{T}\bv{W}(\bv{A}\bv{x}_{k-1}+\bv{B}\bv{u}_{k}-\bv{x}_{d})+0.5\Tr(\bv{\Sigma}_{\bv{x}}\bv{W})$. In turn,  $\E_{\plant{k|k-1}}\left[\hat{c}_{k}(\bv{X}_{k})\right]$ can be conveniently written as: 
\begin{equation}\label{eqn:expected_cost-to-go_Gaussian}
    \begin{split}
      & \E_{\plant{k|k-1}}\left[\hat{c}_{k}(\bv{X}_{k})\right]=0.5\Tr(\bv{S}_{k}\bv{\Sigma})-0.5\omega_{k}\\ 
      & -0.5(\bv{A}\bv{x}_{k-1}-\bv{B}\bv{u}_{k})^{T}\bv{S}_{k}(\bv{A}\bv{x}_{k-1}(\bv{A}\bv{x}_{k-1}-\bv{B}\bv{u}_{k}),
    \end{split}
\end{equation}
with $\omega_{k}$ given, for $k = 1:N-1$, by the following backward recursion (starting with $\omega_{N}=0$)  
\begin{equation*}
    \begin{split}
       &\omega_{k} =\\
          &\ln(|\bv{I}+(\bv{B}\bv{Q}^{0.5})^{T}\bv{\bar{S}}_{k+1}(\bv{B}\bv{Q}^{0.5})|)+\omega_{k+1}+\Tr(\bv{S}_{k+1}\bv{\Sigma})\\
          &+\bv{x}_{d}^{T}(\bv{R}^{-1}+\bv{W})\bv{x}_{d}+\bv{u}_{d}\bv{Q}^{-1}\bv{u}_{d}+\Tr(\bv{R}^{-1}\bv{\Sigma})+\Tr(\bv{W}\bv{\Sigma})\\ 
          &+(\bv{B}^{T}\bv{R}^{-1}\bv{x}_{d}+\bv{Q}^{-1}\bv{u}_{d})^{T}\bv{\Sigma}_{{k}}(\bv{B}^{T}\bv{R}^{-1}\bv{x}_{d}+\bv{Q}^{-1}\bv{u}_{d}), \\ 
           \end{split}
\end{equation*}
and where $\mid\cdot\mid$ denotes the determinant.  In the above expressions ${\bv{S}}_k$ and $\bar{\bv{S}}_k$ are given,  for $k = 1:N-1$,  by the following backward recursion (starting with $\bv{S}_{N} = 0$)
\begin{equation}\label{constraints_cont_var}
    \begin{aligned}
    & \bv{S}_{k} = \bv{A}^{T}\left(\bar{\bv{S}}_{k+1}-\bar{\bv{S}}_{k+1}\bv{B}  \bv{\Sigma}_{k+1}\bv{B}^{T}\bar{\bv{S}}_{k+1}\right)\bv{A},\\
    & \bar{\bv{S}}_{k}=\bv{S}_{k}+\bv{R}^{-1}+\bv{W},\\
     &  \bv{\Sigma}^{\star}_{k} = (\bv{Q}^{-1}+\bv{B}^{T}\bv{\bar{S}}_{k}\bv{B})^{-1}. 
                \end{aligned}
\end{equation}
The policy in \eqref{eqn:gaussian_policy} - \eqref{eqn:gaussian_policy_recursion} is then obtained by substituting in \eqref{eqn:optimal_solution_statement} and \eqref{eqn:p_bar} the above analytical expressions for $\DKL\left(\plant{k|k-1}||\idealplant{k|k-1}\right)$,  $\E_{\plant{k|k-1}}\left[c_{k}(\bv{x}_{k})\right]$ and $\E_{\plant{k|k-1}}\left[\hat{c}_{k}(\bv{X}_{k})\right]$.  Specifically,  adapting the derivations of e.g.,  \cite{KARNY19961719},  by completing the squares we get:
\begin{equation*}
    \begin{aligned}
    \optimalcontrol{k|k-1} & = (2\pi|\bv{\Sigma}_{k})^{-0.5}\exp\big(-0.5\big[(\bv{u}_{k}\\
    &+(\bv{\Sigma}_{k}(\bv{x}_{k-1}^{T}\bv{A}^{T}\bv{\bar{S}}_{k}\bv{B}-\bv{Q}^{-1}\bv{u}_{d}-\bv{x}_{d}\bv{\bar{S}_k}\bv{B})^T\bv{\Sigma}_{k}^{-1}\\
     &(\bv{u}_{k}+(\bv{\Sigma}_{k}(\bv{x}_{k-1}^{T}\bv{A}^{T}\bv{\bar{S}}_{k}\bv{B}-\bv{Q}^{-1}\bv{u}_{d}-\bv{x}_{d}\bv{\bar{S}_k}\bv{B})\big]\big),
    \end{aligned}
\end{equation*}
which is in fact  the Gaussian in \eqref{eqn:gaussian_policy} with $\boldsymbol{\mu}^{\star}_{k}$ and  $\boldsymbol{\Sigma}_k^\star$ given by  \eqref{eqn:gaussian_policy_recursion}. Note that $\omega_{k}$ does not appear in the backward recursion \eqref{eqn:gaussian_policy_recursion} since it  simplifies in the calculations. }

\end{document}